\documentclass[12pt,a4paper]{article}

\usepackage{amssymb}
\usepackage{amsmath}
\usepackage{graphicx}
\def\cM{{\cal M}}
\def\oM{{\overline{\cal M}}}

\def\qed{{\hfill $\diamondsuit$}}

\def\CP{{{\mathbb C}{\rm P}}}
\def\P{{\mathbb P}}

\def\Z{{\mathbb Z}}
\def\C{{\mathbb C}}

\def\E{\mathcal{E}}
\newcommand{\mcO}{\mathcal{O}}
\def\LM{{\rm LM}}
\def\tb{{\tilde b}}
\def\la{\left\langle}
\def\ra{\right\rangle}

\def\DR{{\rm DR}}
\def\DRa{{\rm DR}^{\rm adm}}
\def\DRs{{\rm DR}^{\rm stab}}

\def\disc{{\rm disc}}
\def\ta{{\widetilde a}}
\def\virt{{\rm virt}}

\usepackage{theorem}

\newtheorem{theorem}{Theorem}
\newtheorem{proposition}{Proposition}[section]
\newtheorem{corollary}[proposition]{Corollary}
\newtheorem{lemma}[proposition]{Lemma}
\newtheorem{conjecture}[proposition]{Conjecture}
{\theorembodyfont{\rmfamily}
\newtheorem{definition}[proposition]{Definition}
\newtheorem{example}[proposition]{Example}
\newtheorem{remark}[proposition]{Remark}
\newtheorem{notation}[proposition]{Notation}
}

\title{Integrals of $\psi$-classes over double ramification cycles}
\author{A. Buryak, S. Shadrin, L. Spitz, D. Zvonkine} 

\begin{document}

\maketitle

\begin{abstract}
DR-cycles are certain cycles on the moduli space of curves. Intuitively, they parametrize curves that allow a map to~$\P^1$ with some specified ramification profile over two points. They are known to be tautological classes, but in general there is no known expression in terms of standard tautological classes. In this paper, we compute the intersection numbers of those DR-cycles with any monomials in $\psi$-classes when this intersection is zero-dimensional.
\end{abstract}

\tableofcontents

\section{Introduction}

\subsection{Relative stable maps and double ramification cycles}

Let $a_1, \dots, a_n$ be a list of integers satisfying $\sum a_i = 0$. To a list like that we assign a space of ``rubber'' stable maps to $\CP^1$ relative to $0$ and~$\infty$ in the following way.

Denote by $n_+$ the number of positive integers among the $a_i$'s. They form a partition $\mu = (\mu_1, \dots, \mu_{n_+})$. Similarly, denote by $n_-$ the number of negative integers among the $a_i$'s. After a change of sign they form another partition $\nu = (\nu_1, \dots, \nu_{n_-})$. Both $\mu$ and $\nu$ are partitions of the same integer 
$$
d = \frac12 \sum_{i=1}^n |a_i|.
$$
Finally, let $n_0$ be the number of vanishing $a_i$'s. 

\begin{notation}
To the list $a_1, \dots, a_n$ we assign the space 
$$
\oM_{g;a_1, \dots, a_n} = \oM^{\sim}_{g,n_0; \mu, \nu}(\CP^1, 0, \infty)
$$ 
of degree~$d$ ``rubber'' stable maps to $\CP^1$ relative to $0$ and $\infty$ with ramification profiles $\mu$ and $\nu$, respectively. Here ``rubber'' means that we factor the space by the $\C^*$ action in the target~$\CP^1$. We consider the pre-images of $0$ and $\infty$ as marked points and there are $n_0$ more additional marked points.
\end{notation}

Thus in the source curve there are $n$ numbered marked points with labels $a_1, \dots, a_n$. The relative stable map sends the points with positive labels to $0$, those with negative labels to~$\infty$, while those with zero labels do not have a fixed image.

We have a forgetful map
$$
p \colon \oM_{g;a_1, \dots, a_n} \to \oM_{g,n}.
$$

\begin{definition}
The push-forward 
$$
p_*[ \oM_{g;a_1, \dots, a_n}]^\virt
$$
of the virtual fundamental class under the forgetful map~$p$ is called a {\em double ramification cycle} or a {\em DR-cycle} and is denoted by $\DR_g(a_1, \dots, a_n)$. 
\end{definition}

It is known (see~\cite{FabPan05}) that the Poincar\'e dual cohomology class of $\DR_g(a_1, \dots, a_n)$ lies in the tautological Chow ring of $\oM_{g,n}$. 
The virtual dimension of $\oM_{g;a_1, \dots, a_n}$ and hence the dimension of $\DR_g(a_1, \dots, a_n)$ equals $2g-3+n$. 

\paragraph{Problem.} Find an explicit expression for the class $\DR_g(a_1, \dots, a_n)$.

\bigskip

This is a well-known problem publicized in particular by Y.~Eliashberg in view of applications to Symplectic Field Theory. Recently R.~Hain~\cite{Hain11} found the restriction of $\DR_g(a_1, \dots, a_n)$ to the locus $\oM_{g,n}^{\rm c}$ of curves with compact Jacobians. His expression is a homogeneous polynomial of degree~$2g$ in $a_1, \dots, a_n$ with coefficients in $H^g(\oM_{g,n}^{\rm c})$. In this paper we find the intersection numbers of $\DR_g(a_1, \dots, a_n)$ with monomials in $\psi$-classes. Note that these numbers involve more than the knowledge of $\DR_g(a_1, \dots, a_n)$ on $\oM_{g,n}^{\rm c}$. Thus our results are in some sense complementary with Hain's, even though they are still insufficient to deduce the complete expression for the double ramification cycles. For a given monomial in $\psi_1, \dots, \psi_n$ the intersection number we find is a non-homogeneous polynomial of degree $2g$ in variables $a_1, \dots, a_n$. This gives additional evidence to the following folklore conjecture.

\begin{conjecture}
$\DR_g(a_1, \dots, a_n)$ is a polynomial in $a_1, \dots, a_n$ with coefficients in $H^g(\oM_{g,n})$.
\end{conjecture}

We start with the intersection number of DR-cycles with the power of one $\psi$-class.

\begin{notation} Denote by $S(z)$ the power series
$$
S(z) = \frac{\sinh(z/2)}{z/2} = \sum_{k \geq 0} \frac{z^{2k}}{2^{2k} \, (2k+1)!}   =1 + \frac{z^2}{24} + \frac{z^4}{1920} + \frac{z^6}{322560} + \dots.
$$
\end{notation}

\begin{theorem} \label{Thm:1}
We have
$$
\psi_s^{2g-3+n} \DR_g(a_1, \dots, a_n)
= [z^{2g}]
\frac{\prod_{i \not= s} S(a_i z)}{S(z)},
$$
where $[z^{2g}]$ denotes the coefficient of $z^{2g}$.
\end{theorem}

A generalization of this theorem to DR-cycles with forgotten points is given in Proposition~\ref{Prop:forgottenpoints}. It turns out that for a forgotten marked point the factor $S(a_iz)$ must be replaced by $S(a_iz)-1$.

\subsection{Several $\psi$-classes}

Our next goal is to express the integral over a DR-cycle of a monomial in $\psi$-classes at different marked points. We will use the following notation.

\begin{itemize}
\item
We let $\zeta(z) = e^{z/2} - e^{-z/2}$. (In the previous section we were using $S(z) = \zeta(z)/z$, but here $\zeta(z)$ is much more convenient.)

\item
For a permutation $\sigma \in S_n$ denote $a_i' = a_{\sigma(i)}$ and $z'_i = z_{\sigma(i)}$. 

\item
Finally, 
$$
\begin{array}{|cc|}
a & b
\\ c & d
\end{array}
= ad-bc.
$$
\end{itemize}

\begin{theorem} \label{Thm:2} Given a list of $n$ integers $a_1, \dots, a_n$, satisfying $\sum a_i = 0$ and a list of non-negative integers $d_1, \dots, d_n$ satisfying $\sum d_i = 2g-3+n$, the integral 
$$
\DR_g(a_1, \cdots a_n) \psi_1^{d_1} \cdots \psi_n^{d_n}
$$
of a monomial in $\psi$-classes over a DR-cycle is equal to the coefficient of 
$$
z_1^{d_1} \cdots z_n^{d_n}
$$
in the generating function
\begin{multline*}
\phantom{ab}\frac
{
z_1 \cdots z_n}
{\zeta(z_1 + \cdots + z_n)}
 \sum_{
\substack{\sigma \in S_n\\ \sigma(1) = 1}
}
\\
\frac
{
\zeta\left(
\begin{array}{|cc|}
a_1' & a_2'
\\ z_1' & z_2'
\end{array}
\right)
\, 
\zeta\left(
\begin{array}{|cc|}
a_1'+a_2' & a_3'
\\ z_1'+z_2' & z_3'
\end{array}
\right)
\,
\cdots
\,
\zeta\left(
\begin{array}{|cc|}
a_1' + \cdots + a_{n-1}' & a_n'
\\ z_1' + \cdots + z_{n-1}' & z_n'
\end{array}
\right)
}
{
z_1' 
\;
\begin{array}{|cc|}
a_1' & a_2'
\\ z_1' & z_2'
\end{array}
\; \; 
\begin{array}{|cc|}
a_2' & a_3'
\\ z_2' & z_3'
\end{array}
\;
\cdots
\;
\begin{array}{|cc|}
a_{n-1}' & a_n'
\\ z_{n-1}' & z_n'
\end{array}
\;
z_n'
}.
\end{multline*}
\end{theorem}

\begin{remark}
The expression for the generating function is not written in a symmetrical form: the first marked point is singled out, since we only sum over the permutations that fix the element~1. However the generating function turns out to be symmetric in all $n$ variables. The expression can be symmetrized by extending the summation to all permutations and dividing by~$n$.
\end{remark}

\begin{remark}
At first sight it appears that the generating function has simple poles along the hyperplanes $a_i z_j - a_j z_i$ (because of the determinants in the denominator) and $z_1 + \cdots + z_n = 0$ (because of the $\zeta(z_1 + \cdots + z_n)$ in the denominator). It is easy to see, however, that these denominators actually simplify. 

Indeed, in each summand the factor $a_1'z_2' - a_2' z_1'$ simplifies with $\zeta(a_1'z_2' - a_2' z_1')$. But this was the only factor of the form $a_1 z_i - a_i z_1$, thus no factor like that remains in the denominator of any summand and hence of the total sum. Since the first marked point was singled out arbitrarily, this implies that no factor of the form $a_i z_j - a_j z_i$ remains in the denominator. 

As for the factor $z_1 + \cdots + z_n$, it simplifies with 
$$
\zeta\left(
\begin{array}{|cc|}
a_1' + \cdots + a_{n-1}' & a_n'
\\ z_1' + \cdots + z_{n-1}' & z_n'
\end{array}
\right),
$$
if we take into account that $a_1'+ \cdots + a_{n-1}' = -a_n'$.

The only case where this proof breaks down is when $n=2$. Indeed, in this case $z_1+z_2$ and $a_1 z_2 - a_2 z_1 = a_1 (z_1+z_2)$ are twice the same factor, but this factor is only compensated for once in the numerator. In this case the generating function does contain a singularity of the form $1/(z_1+z_2)$ (see Example~\ref{Ex:nequals2}). This singular term should be ignored when we extract the coefficients. 
\end{remark}

\begin{example} \label{Ex:nequals2}
For $n=2$ we let $a_1 = a$, $a_2 = -a$. There is only one permutation in $S_2$ that fixes the first element. Thus we get the generating function
$$
\frac{z_1 z_2}{\zeta(z_1+z_2)} \; \frac{\zeta(a(z_1+z_2))}{z_1 \, a(z_1+z_2) \, z_2} = \frac{\zeta(a(z_1+z_2))}{a (z_1 + z_2) \zeta(z_1+z_2)}
$$
$$
= \frac1{z_1+z_2} + \frac{a^2-1}{24} (z_1 + z_2) + \frac{(a^2-1)(3a^2-7)}{5760} (z_1+z_2)^3+ \cdots
$$
It follows that
\begin{align*}
\DR_1(a,-a) \psi_1 &= \DR_1(a,-a) \psi_2 = \frac{a^2-1}{24},\\
\DR_2(a,-a) \psi_1^3 &= \DR_2(a,-a) \psi_2^3 = \frac{(a^2-1)(3a^2-7)}{5760},\\
\DR_2(a,-a) \psi_1^2 \psi_2 &= \DR_2(a,-a) \psi_1 \psi_2^2 \\
&= \frac{3(a^2-1)(3a^2-7)}{5760} = \frac{(a^2-1)(3a^2-7)}{1920}.
\end{align*}
\end{example}

\begin{example} \label{Ex:nequals3}
For $n=3$ we have $a_3 = -(a_1+a_2)$. There are two summands in the formula corresponding to the permutations $(1,2,3)$ and $(1,3,2)$. We get
$$
\frac1{\zeta(z_1 + z_2 + z_3)} 
\left\{
\frac{\zeta(a_1z_2 - a_2 z_1)}{a_1z_2 - a_2 z_1} 
\frac{z_2\, \zeta\bigl((a_1+a_2)(z_1+z_2+z_3)\bigr)}
{a_2z_3+(a_1+a_2)z_2}
+
\right.
$$
$$
\left.
\frac{\zeta(a_1z_3 + (a_1+a_2) z_1)}{a_1z_3 +(a_1+a_2) z_1} 
\frac{z_3\, \zeta\bigl(a_2(z_1+z_2+z_3)\bigr)}
{a_2z_3+(a_1+a_2)z_2}
\right\}
.
$$
Expanding this expression we get, in particular,
$$
\DR_1(a_1, a_2, a_3) \psi_1^2 = \frac{a_2^2+a_3^2-1}{24},
$$
$$
\DR_1(a_1, a_2, a_3) \psi_1 \psi_2 = \frac{a_1^2+a_2^2+a_3^2-2}{24},
$$
where we have re-introduced $a_3$ for more symmetry.

\end{example}

\subsection{Completed cycles as a particular case of Theorem~\ref{Thm:1}}
\label{Ssec:CompCyc}

Let $\oM^0_{g,1,K;\kappa}(\CP^1, \infty)$ be the space of degree~$K$ relative stable maps $f \colon C\to \CP^1$ with branching profile $\kappa = (k_1, \dots, k_n)$ over~$\infty$ and with one marked point $x \in C$ satisfying the condition that $f(x) = 0$. It is a natural problem to find an effective cycle representing the homology class
$$
[\oM^0_{g,1,K; \kappa}(\CP^1, \infty)]^\virt \; \psi_x^m.
$$
Okounkov and Pandharipande gave an answer to this question when $m$ is equal to the virtual dimension $K+n+2g-2$ of $\oM^0_{g,1,K;\kappa}(\CP^1,\infty)$ and thus the answer is just a number. To simplify the formula we assume that the $n$ pre-images of $\infty$ in our space of relative stable maps are numbered. Then we have the following equality.

\paragraph{Theorem} (Okounkov, Pandharipande~\cite{OkoPan06}). 
For $m = K+n+2g-2$ we have
$$
[\oM^0_{g,1,K; \kappa}(\CP^1, \infty)]^\virt \; \psi_x^m = 
m! \frac{\prod_{i=1}^n k_i}{K!} [z^{2g}] S(z)^{K-1} \prod_{i=1}^n S(k_i z).
$$

Using the degeneration of the target it is not hard to generalize this formula to several relative points with ramification types $\mu_1, \dots, \mu_s$. In particular, for the case of two relative points, the following expression is given in~\cite{OkoPan06}, Eq.~(3.11) or \cite{Rossi}, Eq.~(10). Let $a_1, \dots, a_n$ be the list of elements of $\mu_1$ merged with the list of elements of $\mu_2$ with reversed signs. Thus $\sum a_i = 0$. Denote by $\psi_x$ the $\psi$-class at the marked point~$x$.

\paragraph{Theorem} (Okounkov-Pandharipande, Rossi) We have
$$
[\oM^0_{g,1,d; \mu_1, \mu_2}(\CP^1, p_1, p_2)]^\virt \; \psi_x^{n+2g-2} 
= [z^{2g}]
\frac{\prod_{i=1}^n S(a_i z)}{S(z)}.
$$

It is easy to see that this formula is a particular case of Theorem~\ref{Thm:1}, namely, the case when $a_s = 0$ while all other $a_i$'s do not vanish. The case where $a_s=0$ and some other $a_i$'s may also vanish is covered by a more general computation in Proposition~2.5 of~\cite{OkoPan06b}.

Actually, we don't have an independent proof for the case $a_s=0$ we just invoke the above result. Our proof for the case $a_s \not=0$ is quite different and doesn't generalize to $a_s=0$. Thus we get the same answer for $a_s=0$ and $a_s \not=0$, even though we do not know any proof that would work in both situations.

\subsection{A conjecture on $\psi^m$}

Here we show that our formulas are compatible with a conjectural identity on virtual fundamental classes. This section can be skipped in first reading.

Let us return to the opening question of Section~\ref{Ssec:CompCyc}: find an explicit cycle that represents $\psi_x^m$ in the space of stable maps. For simplicity we consider here the space~$\oM^0_{G,1,d}(\CP^1)$ of genus~$G$ degree~$d$ stable maps~$f$ with no relative points and just one marked point~$x$ satisfying $f(x) = 0$.

\begin{notation}
Let $\DR_g(a_1, \dots, a_n)$ be a DR-cycle in $\oM_{g,n}$ and $\pi \colon \oM_{g,n} \to \oM_{g,n-1}$ the forgetful map forgetting, for instance, the $n^{\mathrm{th}}$ point. We will denote the image of the DR-cycle in $\oM_{g,n-1}$ by
$\DR_g(a_1, \dots, a_{n-1}; \ta_n)$. Similarly, if we take the image of a DR-cycle under a map forgetting several marked points we will put tildes over the corresponding $a_i$'s.
\end{notation}

To formulate the conjecture, let us first describe a particular homology class in~$\oM^0_{G,1,d}(\CP^1)$ shown in the figure.
\begin{center}
\ 
\includegraphics[width=20em]{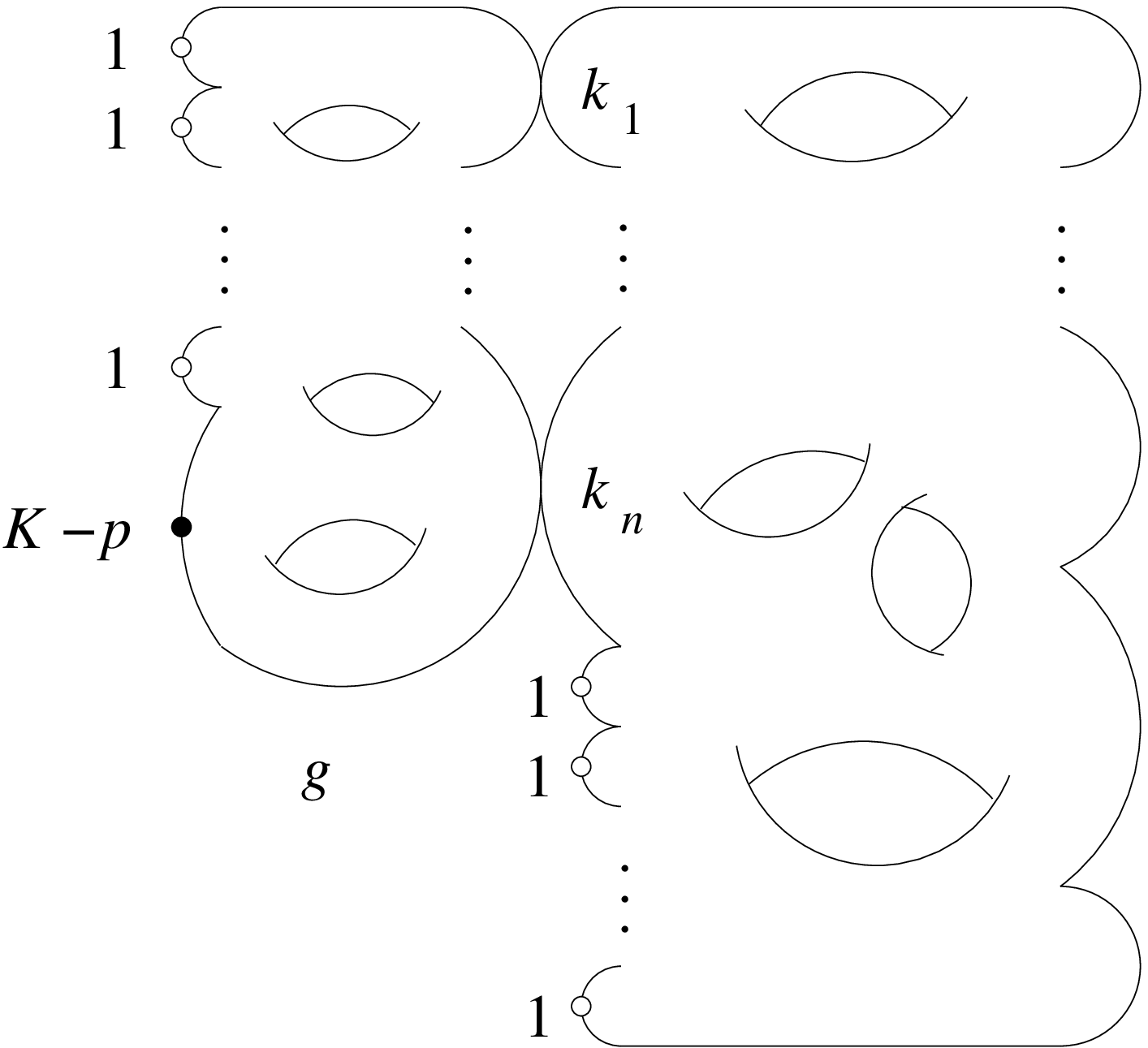}

\end{center}

The left-hand side of the figure represents the space of rubber maps of genus~$g$ with branching profiles $(k_1, \dots, k_n)$ and $(K-p, 1, \dots, 1)$. 
The right-hand side of the figure represents the space of stable maps to $\CP^1$ relative to $0$. The branching profile over~$0$ is given by the partition $(k_1, \dots, k_n, 1, \dots, 1)$. The source curve in the right-hand side is not necessarily connected and has Euler characteristic $2-2G'$, where~$G' = G-g-n+1$. Multiplying the virtual fundamental classes of both spaces we obtain a virtual fundamental class on the product space.

Now, there is a natural forgetful map from the product space to the space~$\oM^0_{G,1,d}(\CP^1)$. It is obtained by contracting all the rubber components in the target rational curve, forgetting the marked points shown in white on the picture, and stabilizing the resulting map. We will denote the image of the virtual fundamental class under this map by
$$
\DR_g\left(\framebox{K-p}, -k_1, \dots, -k_n; \; {\tilde 1}^p \right) 
\; \boxtimes \;
\left[\oM^\disc_{G',0,d; \kappa}(\CP^1,0) \right]^\virt.
$$
Here ${\tilde 1}^p$ means that there are $p$ points marked with the integer $1$ that are forgotten by the forgetful map; the boxed number corresponds to the only marked point that survives after we apply the forgetful map; $\kappa$ is the partition $(k_1, \dots, k_n, 1^{d-K})$. Even though the notation is pretty complicated, it is actually still incomplete: indeed, in the second factor what we mean is not the virtual fundamental class of the space itself, but its image under the forgetful map to the space~$\oM^0_{G,1,d}(\CP^1)$.

\begin{conjecture} \label{Conj:psitothem}
We have
\begin{align*}
\lefteqn{m! \psi_x^m \; [\oM^0_{G,1,d}(\CP^1)]^\virt  =} 
\qquad \\
& 
\lefteqn{\sum_{n \geq 1}
\frac1{n!}
\sum_{g=0}^{G-n+1}
\sum_{p=0}^g
\sum_{K=p}^{m+p}
\frac{m!}{p!\,(K-p)!} \, \psi_x^{m-K+p}
\sum_{
\substack{k_1, \dots, k_n\\
\sum k_i = K}
}
\prod_{i=1}^n k_i \times}
\qquad
\\
&&  \DR_g\left(\framebox{K-p}, -k_1, \dots, -k_n; \; {\tilde 1}^p\right) 
\; \boxtimes \;
\left[\oM^\disc_{G',0,d; \kappa}(\CP^1,0) \right]^\virt,
\end{align*}
where $G' = G-g-n+1$, $\kappa = (k_1, \dots, k_n, 1^{d-K})$.
\end{conjecture}

This conjecture is proved in genus~0 (see~\cite{KaLaZv}). It is also known that the conjecture is true for any~$m$ if it is true for $m=0$ (unpublished). For $m=0$ the conjecture describes a splitting of the virtual fundamental class of $\oM_{G,1,d}^0$ according to the type of singularity of the map $f$ at the marked point~$x$. 

Now we add another piece of evidence in favour of the conjecture.

\begin{theorem} \label{Thm:evidence}
Let $g \geq 0$, $n \geq 1$, and $k_1, \dots, k_n$. Denote $K = \sum k_i$ and assume let $m = K+n+2g-2$. Then we have
$$
\sum_{p=0}^g
\frac{m!}{p! \, (K-p)!} \, \psi_x^{m-K+p} \;
\prod_{i=1}^n k_i\; 
\DR_g\left(\framebox{K-p}, -k_1, \dots, -k_n; \; {\tilde 1}^p\right) = 
$$
$$
m! \frac{\prod_{i=1}^n k_i}{K!} [z^{2g}] S(z)^{K-1} \prod_{i=1}^n S(k_i z),
$$
where the integer corresponding to the marked point~$x$ is boxed.
\end{theorem}

The left-hand side of the equality extracts the numerical terms from the conjectural expression of $\psi^m$, that is, those terms where the power of the $\psi$-class matches the dimension of the DR-cycle. The right-hand side of the equality coincides with Okounkov and Pandharipande's expression for the coefficients of the completed cycles (see Section~\ref{Ssec:CompCyc}). Thus the theorem shows that the conjecture allows one to recover Okounkov and Pandharipande's result enriching it with non-numerical terms that are invisible in their formulas.

\subsection{Acknowledgements}
The authors would like to thank J.~Wise, R.~Pandharipande and J.~Li for useful discussions. A.~B, S.~S and L.~S were supported by a Vidi grant of the Netherlands Organization for Scientific Research. A.~B was further supported by Russian Federation Government grant no. 2010-220-01-077 (ag. no. 11.634.31.0005), the grants RFBR-10-010-00678, NSh-4850.2012.1 and the Moebius Contest Foundation for Young Scientist. D.~Z was supported by the grant ANR-09-JCJC-0104-01.

\section{DR-cycle times a $\psi$-class}

\subsection{The splitting formulas - formulations}

In this section we express the product of $\DR_g(a_1, \dots, a_n)$ and the class $\psi_s$ for some $s\in \{ 1, \dots, n \}$ in terms of other DR-cycles. This will make it possible to evaluate monomials in $\psi$-classes on a DR-cycle by induction. Note that we can only do that if $a_s \ne 0$. Conjecture~\ref{Conj:psitothem} is intended to be a generalization of our recursive formulas to the case $a_s=0$, but it is still open.

The picture below shows a cycle in $\oM_{g,n}$ obtained from two DR-cycles via a gluing map. 

The two DR-cycles are constructed in the following way. The list $a_1, \dots, a_n$ is divided into two disjoint parts: $I \sqcup J = \{1, \dots, n \}$ in such a way that $\sum_{i \in I} a_i >0$ or, equivalently, $\sum_{i \in J} a_i <0$. In the figure, for instance, we have $1, 2, s \in I$ and $n \in J$. Then a new list of positive integers $k_1, \dots, k_p$ is chosen in such a way that 
$$
\sum_{i \in I} a_i - \sum_{i=1}^p k_i = \sum_{i \in J} a_i + \sum_{i=1}^p k_i = 0.
$$
Now two DR-cycles of genera $g_1$ and $g_2$ are formed as shown in the figure and glued together at the ``new'' marked points labelled $k_1, \dots, k_p$. Since we want to get a genus~$g$ in the end we impose the condition $g_1 + g_2 + p-1 = g$. We denote by 
$$
\DR_{g_1}(a_I,-k_1, \dots, -k_p) \boxtimes \DR_{g_2}(a_J, k_1, \dots, k_p)
$$
the resulting cycle in~$\oM_{g,n}$.

\begin{center}
\ 
\includegraphics[width=15em]{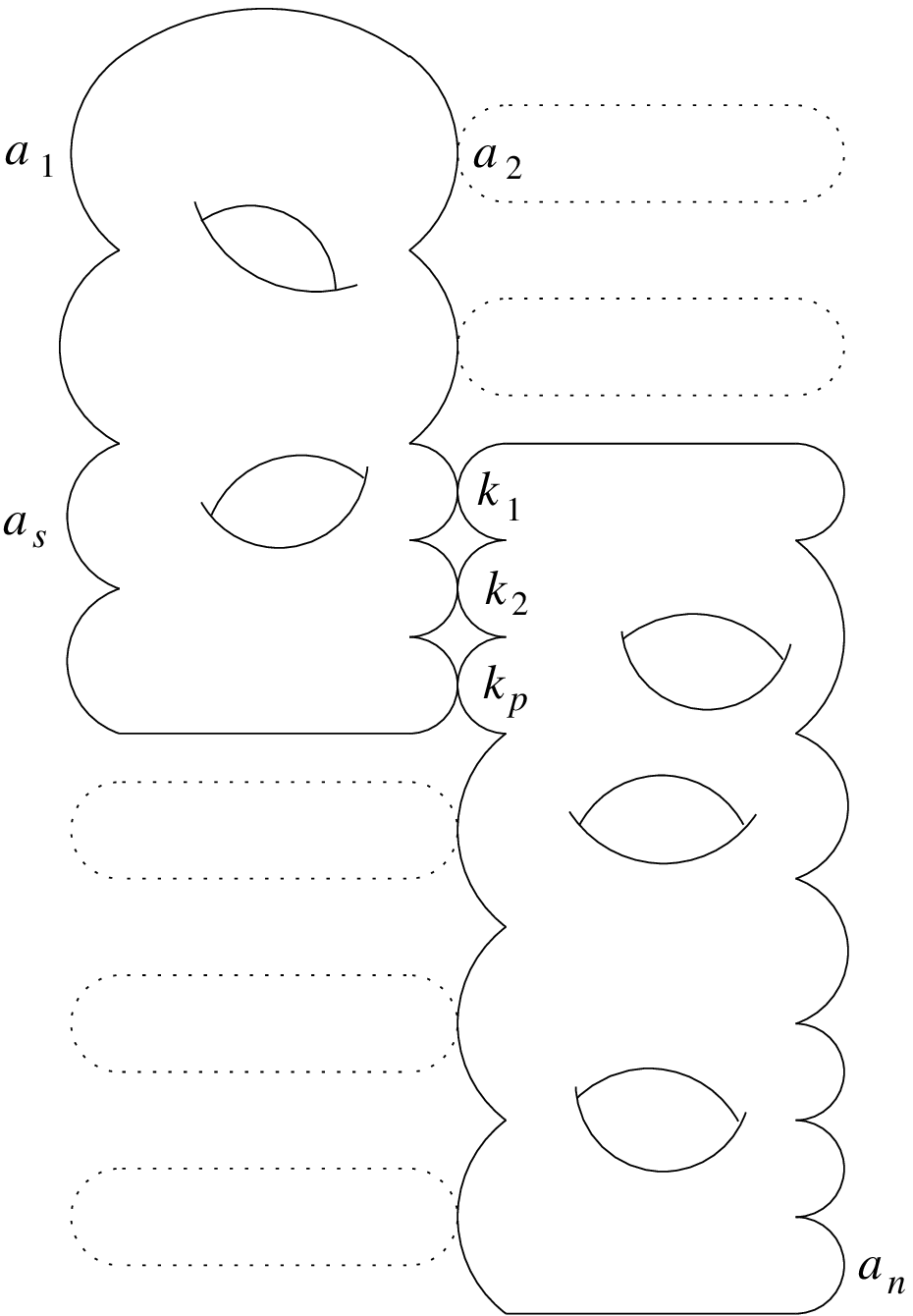}
\end{center}

Let $r = 2g-2+n$ be the number of branch points of our initial DR-cycle $\DR(a_1, \dots, a_n)$. Let $r' = 2g_1 - 2 + |I|+p$ and $r'' = 2g_2-2+|J|+p$ be the numbers of branch points in the two components of the target curve. (In both cases we don't count $0$ and $\infty$.)

\begin{theorem} \label{Thm:geomlemma}
Let $a_1, \dots, a_n$ be a list of integers with vanishing sum. Assume that $a_s \not= 0$. Then we have
$$
a_s \psi_s \DR_g(a_1, \dots, a_n) = 
$$
$$
\sum_{I,J}
\sum_{p \geq 1}\sum_{g_1, g_2} 
\sum_{k_1, \dots, k_p}
 \frac{\rho}{r} \frac{\prod_{i=1}^p k_i}{p!}
\DR_{g_1}(a_I,-k_1, \dots, -k_p) \boxtimes \DR_{g_2}(a_J, k_1, \dots, k_p).
$$
Here the first sum is taken over all $I \sqcup J = \{ 1, \dots, n\}$ such that $\sum_{i \in I} a_i >0$; the third sum is over all non-negative genera $g_1$, $g_2$ satisfying
$g_1 + g_2 + p-1= g$; the fourth sum is over the $p$-uplets of positive integers with total sum $\sum_{i \in I} a_i = -\sum_{i \in J} a_i$. The number $\rho$ is defined by
$$
\rho = \left\{ 
\begin{array}{lcl}
 r'' & \mbox{if} & s \in I, \\
-r' & \mbox{if} & s \in J.
\end{array}
\right.
$$
\end{theorem}

\begin{theorem} \label{Thm:geomlemmabis}
Let $a_1, \dots, a_n$ be a list of integers with vanishing sum. Assume that $a_s \not= 0$ and $a_l=0$. Then we have
$$
a_s \psi_s \DR_g(a_1, \dots, a_n) = 
$$
$$
\sum_{I,J}
\sum_{p \geq 1}\sum_{g_1, g_2} 
\sum_{k_1, \dots, k_p}
\varepsilon \; \frac{\prod_{i=1}^p k_i}{p!}
\DR_{g_1}(a_I,-k_1, \dots, -k_p) \boxtimes \DR_{g_2}(a_J, k_1, \dots, k_p).
$$
Here the first sum is taken over all $I \sqcup J = \{ 1, \dots, n\}$ such that $\sum_{i \in I} a_i >0$; the third sum is over all non-negative genera $g_1$, $g_2$ satisfying
$g_1 + g_2 + p-1= g$; the fourth sum is over the $p$-uplets of positive integers with total sum $\sum_{i \in I} a_i = -\sum_{i \in J} a_i$. The number $\varepsilon$ is defined by
$$
\varepsilon = \left\{ 
\begin{array}{rl}
 1 & \mbox{if} \quad s \in I, l \in J,\\
-1 & \mbox{if} \quad  s \in J, l \in I,\\
0 & \mbox{otherwise.}
\end{array}
\right.
$$
\end{theorem}

We will call these Theorem~\ref{Thm:geomlemma} the {\em splitting formula with respect to branching points} and Theorem~\ref{Thm:geomlemmabis} the {\em splitting formula with respect to a marked point}.  Before proving the theorems let us formulate some corollaries that we will use in our computations.

\begin{corollary} \label{Cor:onepsi}
Assume that $a_s \not=0$. We have
\begin{align*}
&r a_s \; \psi_s^{r-1} \DR_g(a_1, \dots, a_n)
= \\
&- \frac12 \sum_{i,j \not= s}
(a_i+a_j) \; \psi_s^{r-2} \DR_g(a_1, \dots, \widehat{a_i}, \dots, \widehat{a_j}, \dots, a_n, a_i+a_j)\\
&-\frac12 \sum_{i \not= s} 
{\rm sign}(a_i) \sum_{\substack{b+c=a_i\\ b \cdot c>0}}
bc \; \psi_s^{r-2} \DR_{g-1}(a_1, \dots, \widehat{a_i}, \dots, a_n, b, c).
\end{align*}
Here, as before, $r=2g-2+n$ and a hat means that the element is skipped. 
\end{corollary}

\paragraph{Proof.} We will use the splitting formula with respect to the branch points. Since we are interested in the intersection number of our DR-cycle with $\psi_s^{r-1}$ we only need to keep those terms of the splitting formula for which the $s$th marked point stays on a DR-cycle of dimension~$r-2$. This implies that the remaining DR-cycle is of dimension~$0$, that is, it is of the form $\DR_0(a,b,c)$. The expression in the corollary is a sum over all splittings of this form. \qed

\bigskip

This corollary gives a recursive relation for intersection numbers of DR-cycles with powers of {\em one} $\psi$-class. We will use it to prove Theorem~\ref{Thm:1}.

\begin{corollary} \label{Cor:movepsi}
Let $t$ and $s$ be two different elements in $\{ 1, \dots, n \}$. Assume that both $a_s$ and $a_t$ are non-zero. Then we have
$$
(a_s \psi_s - a_t \psi_t) \DR_g(a_1, \dots, a_n)
$$
$$
= \sum_{s \in I, t \in J}
\sum_{p \geq 1}\sum_{g_1, g_2} 
\sum_{k_1, \dots, k_p}
\frac{\prod_{i=1}^p k_i}{p!}
\DR_{g_1}(a_I,-k_1, \dots, -k_p) \boxtimes \DR_{g_2}(a_J, k_1, \dots, k_p)\\
$$
$$ 
-
\sum_{t \in I, s \in J}
\sum_{p \geq 1}\sum_{g_1, g_2} 
\sum_{k_1, \dots, k_p}
\frac{\prod_{i=1}^p k_i}{p!}
\DR_{g_1}(a_I,-k_1, \dots, -k_p) \boxtimes \DR_{g_2}(a_J, k_1, \dots, k_p).
$$
Here, as before, the first sum is taken over all $I \sqcup J = \{ 1, \dots, n\}$ such that $\sum_{i \in I} a_i >0$; the third sum is over all non-negative genera $g_1$, $g_2$ satisfying
$g_1 + g_2 + p-1= g$; the fourth sum is over the $p$-uplets of positive integers with total sum $\sum_{i \in I} a_i = -\sum_{i \in J} a_i$.
\end{corollary}

\paragraph{Proof.} This follows directly from the splitting formula with respect to the branch points. It suffices to notice that the expressions it provides for $a_s \psi_s$ and $a_t \psi_t$ only differ in the definition of $r'$. \qed

\bigskip

Multiplying the identity in this corollary by any monomial in $\psi$-classes of degree $2g-4+n$ we obtain a simple way to ``move'' a $\psi$-class from one marked point to another.

\subsection{The splitting formulas - proofs}

\subsubsection{Plan of proof}
Our proof uses the Losev-Manin compactification $\LM_r$ of $\cM_{0,r+2}$. It is the moduli space of chains of spheres with two special ``white'' marked points $0$ and $\infty$ at the extremities of the chain and $r$ more ``black'' marked points on the other spheres. The black points are allowed to coincide with each other and there should be at least one black point per sphere. For more details see~\cite{LosMan}.
\begin{center}
\ 
\includegraphics[width=15em]{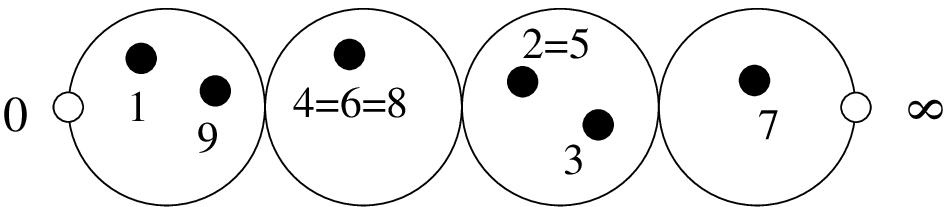}
\end{center}
We have two forgetful maps from the DR-space $\oM_{g;a_1, \dots, a_n}$:
$$
\LM_{r+n_0}/S_r \stackrel{q}{\longleftarrow} \oM_{g;a_1, \dots, a_n} \stackrel{p}{\longrightarrow} \oM_{g,n},
$$
where $n_0$ is the number of indices $i$ such that $a_i=0$ and $r = 2g-2+n$ is the number of branch points.

The map $q$ assigns to a relative stable map its target rational curve. The marked points are the $r$ branch points and the images of the marked points in the source curve. The map~$p$ assigns to a relative stable map its stabilized source curve. (This is the map that we used to define the DR-cycle $\DR_g(a_1, \dots, a_n)$.)

The proof of the splitting formulas proceeds as follows.
\begin{enumerate}
\item Identify the $\psi$-class on the DR-space with the $\psi$-class on the Losev-Manin space.
\item Express the $\psi$-class on the Losev-Manin space as a sum of boundary divisors (we will do that in two ways, whence two splitting formulas).
\item Lift these divisors to the DR-space
\item Subtract the difference between the $\psi$-class on the DR-space and the $\psi$-class on $\oM_{g,n}$.
\end{enumerate}

We start with two lemmas that will be needed in the course of the proof.

\subsubsection{DR-cycles with disconnected domains}
Consider the space of stable maps to $\CP^1$ relative to $0$ and $\infty$, but with disconnected domains
$$
\oM_{g_1; a^1_1, \dots, a^1_{n_1}} \times \cdots \times \oM_{g_k; a^k_1, \dots, a^k_{n_k}}.
$$
We assume that $2g_i-2+n_i > 0$ for each~$i$. From the corresponding rubber space
$$
\left(\oM_{g_1; a^1_1, \dots, a^1_{n_1}} \times \cdots \times \oM_{g_k; a^k_1, \dots, a^k_{n_k}}\right)^{\sim}
$$
there is a natural forgetful map~$p$ to the product of moduli spaces $\oM_{g_1, n_1} \times \cdots \times \oM_{g_k,n_k}$.

\begin{lemma} \label{Lem:onlyone}
The image of the virtual fundamental class of 
$$
\left(\oM_{g_1; a^1_1, \dots, a^1_{n_1}} \times \cdots \times \oM_{g_k; a^k_1, \dots, a^k_{n_k}}\right)^{\sim}
$$
in $\oM_{g_1, n_1} \times \cdots \times \oM_{g_k,n_k}$ under the forgetful map~$p$ vanishes.
\end{lemma}

Even though the computations of this section take place in the DR-space, the goal of the paper is to study the DR-cycles, that is, the images of the virtual fundamental classes of DR-spaces by the map $p$. Therefore in the sequel of this section we will perform all our computations ``modulo terms with disconnected domains''. In other words, we will disregard all the terms that, according to the lemma, vanish after the push-forward by~$p$.

\paragraph{Proof.} We will call {\em parts} the $k$ connected components of the curves.

\paragraph{Adding a new marked point.}
Consider the space
$$
\left(\oM_{g_1; a^1_1, \dots, a^1_{n_1},0} \times \cdots \times \oM_{g_k; a^k_1, \dots, a^k_{n_k}}\right)^{\sim}.
$$
If $\pi$ is the forgetful map that forgets the new point, we have, by the dilaton relation,
$$
\pi_*\left\{\left[\left(\oM_{g_1; a^1_1, \dots, a^1_{n_1},0} \times \cdots \times \oM_{g_k; a^k_1, \dots, a^k_{n_k}}\right)^{\sim} \; \right]^\virt \psi_{n_1+1} \right\}
$$
$$
= (2g_1 - 2 + n_1) 
\left[\left(\oM_{g_1; a^1_1, \dots, a^1_{n_1}} \times \cdots \times \oM_{g_k; a^k_1, \dots, a^k_{n_k}}\right)^{\sim}\; \right]^\virt.
$$
Thus it suffices to prove that the image of 
$$
\left[\left(\oM_{g_1; a^1_1, \dots, a^1_{n_1},0} \times \cdots \times \oM_{g_k; a^k_1, \dots, a^k_{n_k}}\right)^{\sim} \; \right]^\virt \psi_{n_1+1}
$$
vanishes in  $\oM_{g_1, n_1+1} \times \cdots \times \oM_{g_k,n_k}$.

\paragraph{Introducing a $\C^*$-action}
On the space 
$$
\left(\oM_{g_1; a^1_1, \dots, a^1_{n_1},0} \times \cdots \times \oM_{g_k; a^k_1, \dots, a^k_{n_k}}\right)^{\sim}
$$
we can introduce a $\C^*$-action in the following way. Let $f:C \to S$ be a rubber map, where $S$ is a genus~0 curve from the Losev-Manin space. Let $S_\bullet$ be the irreducible component of~$S$ that contains the image of the new marked point, that is, the $(n_1+1)^{\mathrm{st}}$ marked point in the first part of~$C$. (The purpose of adding a new marked point was precisely to be able to single out a component of~$S$ in this way.) Now, for $\lambda \in \C^*$, we let $\lambda.f$ be equal to $f$ on every component of $C$ that does not map to~$S_\bullet$ or is in the first part (that is, the part that contains the new marked point). On the components of the other parts that map to $S_\bullet$ we let $\lambda.f = \lambda f$. 

The pull-back of any differential form from $\oM_{g_1, n_1+1} \times \cdots \times \oM_{g_k,n_k}$ to our DR-space is $\C^*$-invariant, because the action of $\C^*$ does not change the complex structure of the source curve. We are going to prove by localization that the integral against
$$
\left[\left(\oM_{g_1; a^1_1, \dots, a^1_{n_1},0} \times \cdots \times \oM_{g_k; a^k_1, \dots, a^k_{n_k}}\right)^{\sim} \; \right]^\virt \psi_{n_1+1}
$$
of any $\C^*$-invariant form vanishes.

\paragraph{Localization.}
The invariant locus of the $\C^*$-action is composed of maps that have no marked or ramification points over $S_\bullet$ on parts $2, \dots, k$. Thus the invariant locus has three types of components, classified by the topological type of the target genus~0 curve~$S$ (at the generic point of the component of the locus):
\begin{enumerate}
\item the curve~$S$ has the form $S' \cup S_\bullet$;
\item the curve~$S$ has the form $S_\bullet \cup S''$;
\item the curve~$S$ has the form $S' \cup S_\bullet \cup S''$.
\end{enumerate}
Each component of the invariant locus is the product of two (in the first two cases) or three (in the last case) disconnected DR-spaces and has the same virtual fundamental class. A simple dimension count shows that the virtual dimension of each component of the invariant locus is less than the virtual dimension of the original DR-space. (Indeed, the dimension is equal to the number of marked and branch points minus the number of components of~$S$.) Therefore each term in the localization formula vanishes. 

This proves the lemma. \qed

\subsubsection{Pull-backs of divisors from the Losev-Manin space}

Consider a DR-space $\oM_{g;a_1, \dots, a_n}$ and consider the forgetful map
$$
q: \oM_{g;a_1, \dots, a_n} \to \LM_{r + n_0} / S_r.
$$
Let $\alpha \sqcup \beta$ be a partition of the set of indices $i$ such that $a_i=0$. Let $r'+r'' = r$. Denote by $D_{(r',\alpha | r'', \beta)}$ the boundary divisor in the space $\LM_{r + n_0} / S_r$ with self-explanatory notation.

\begin{lemma} \label{Lem:pullback}
Modulo terms with disconnected domains, we have
$$
q^*D_{(r',\alpha | r'', \beta)}
\left[ \oM_{g;a_1, \dots, a_n} \right]^\virt
=
$$
$$
\sum_{I,J}
\sum_{p \geq 1}\sum_{g_1, g_2} 
\sum_{k_1, \dots, k_p}
\frac{\prod_{i=1}^p k_i}{p!}
\left[\oM_{g_1;a_I,-k_1, \dots, -k_p}\right]^\virt 
\boxtimes 
\left[ \oM_{g_2; a_J, k_1, \dots, k_p} \right]^\virt.
$$
Here the first sum is taken over all $I \sqcup J = \{ 1, \dots, n\}$ such that $\alpha \subset I$, $\beta \subset J$ and $\sum_{i \in I} a_i >0$; the third sum is over all non-negative genera $g_1$, $g_2$ satisfying
$g_1 + g_2 + p-1= g$; the fourth sum is over the $p$-uplets of positive integers with total sum $\sum_{i \in I} a_i = -\sum_{i \in J} a_i$.
\end{lemma}

\paragraph{Proof.} This lemma is a version of Jun Li's degeneration formula~\cite{Li1,Li2}. It should be applied to the target rational curve where all branch points have been marked and numbered. We then take the sum of contributions from all possible ways to put $r'$ marked point on one component of the degeneration and $r''$ on the other component. After this we can forget the numbering of the branch points once again.

In the degeneration formula we see DR-spaces with both connected and disconnected domains over each component of the target. However, since we are working modulo terms with disconnected domains, only the terms indicated in the lemma survive. \qed

\subsubsection{Comparing the $\psi$-classes on different spaces}

Recall the two forgetful maps from the DR-space
$$
p: \oM_{g;a_1, \dots, a_n} \to \oM_{g,n}
$$
and
$$
q: \oM_{g;a_1, \dots, a_n} \to \LM_{r+n_0} / S_r.
$$ 

Denote by $\Psi_s$ and $\psi_s$ the $\psi$-classes at the $s$th marked point on $\oM_{g;a_1, \dots, a_n}$ and on $\oM_{g,n}$ respectively. Denote by $\psi_0$ and $\psi_\infty$ the $\psi$-classes on the Losev-Manin space. 

\begin{proposition} \label{Prop:Ionel}
Assume that $a_s \ne 0$. We have
$$
a_s \Psi_s = q^*\psi_0 \quad \mbox{if} \quad a_s > 0,
$$
$$
-a_s \Psi_s = q^*\psi_\infty \quad \mbox{if} \quad a_s < 0.
$$
\end{proposition}

This simple but very useful statement first appeared in Ionel's paper~\cite{Ion02}. Assume for definiteness that $a_s > 0$. Then $q$ obviously identifies the tangent line to $0$ in the target with the $a_s$th power of the tangent line to the $s$th marked point in the source, which proves the proposition.

\begin{lemma} \label{Lem:Psiminuspsi}
Assume that $a_s \ne 0$. Modulo terms with disconnected domains we have
$$
\Psi_s - p^* \psi_s = 
$$
$$
\frac1{|a_s|}
\sum_{I,J}
\sum_{p \geq 1}\sum_{g_1, g_2} 
\sum_{k_1, \dots, k_p}
\frac{\prod_{i=1}^p k_i}{p!}
\left[\oM_{g_1;a_I,-k_1, \dots, -k_p}\right]^\virt 
\boxtimes 
\left[ \oM_{g_2; a_J, k_1, \dots, k_p} \right]^\virt.
$$
Here the first sum is taken over all $I \sqcup J = \{ 1, \dots, n\}$ such that $\sum_{i \in I} a_i >0$ and $s \in J$ if $a_s >0$ or $s \in I$ if $a_s <0$; the third sum is over all non-negative genera $g_1$, $g_2$ satisfying
$g_1 + g_2 + p-1= g$; the fourth sum is over the $p$-uplets of positive integers with total sum $\sum_{i \in I} a_i = -\sum_{i \in J} a_i$.
\end{lemma}

\paragraph{Proof.} The sum in the lemma enumerates all the boundary divisors, modulo the ones with disconnected domains, on which the $s$th marked point lies on a bubble (that is, on a rational component that gets contracted by the forgetful map~$p$). It is precisely those divisors that contribute to the difference between the two $\psi$-classes. It remains to determine the coefficients. 

Assume, for definiteness, that $a_s > 0$. A divisor enumerated in the sum splits the marked and branch points in the target into two groups: those that lie on the component of 0 and those that lie on the component of~$\infty$. Consider the map~$\tilde{q}$ that forgets all the marked and branch points from the component of~$0$. It is a forgetful map between two Losev-Manin spaces. Denote by $\psi_0'$ the $\psi$-class at $0$ on the smaller Losev-Manin space, that is, on the image of the forgetful map. It is easy to see that in the neighbourhood of our divisor and outside of the other divisors enumerated in the lemma we have
$a_s \psi_s = q^* \tilde{q}^* \psi'_0$ and therefore 
$$
a_s(\Psi_s - p^*\psi_s) = q^*(\psi_0 -\tilde{q}^* \psi'_0).
$$
The difference $\psi_0 - \tilde{q}^* \psi_0'$ is exactly given by the divisor where the target curve degenerates. Therefore the coefficient of our divisor in $\Psi_s - p^* \psi_s$ is equal to the coefficient of the same divisor in Lemma~\ref{Lem:pullback} divided by $a_s$. \qed

\subsubsection{Expressing $\psi_0$ and $\psi_\infty$ as boundary divisors}

The class $\psi_0$ on the Losev-Manin space $\LM_{r+n_0}$ is easily expressed as a sum of boundary divisors: namely, for any $i \in \{1, \dots, r+n_0 \}$, we have
\begin{equation} \label{Eq:psimarked}
\psi_0 = \sum \left[ 
\begin{picture}(7,2)
\put(0,-1){\includegraphics[width=7em]{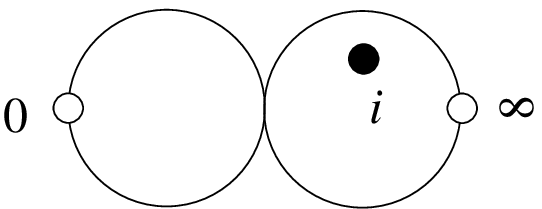}}
\end{picture}
\right], 
\end{equation}
where the sum is over all boundary divisors such that the $i$th marked point lies on the same component as~$\infty$. If $i$ is the image of a marked point we leave this expression as it is.

If $i$ is a branch point, it makes sense to symmetrize the expression with respect to the $S_r$ action, since we are working with the quotient $\LM_{r+n_0} / S_r$. We get 
\begin{equation} \label{Eq:psibranch}
\psi_0 = \sum \frac{r''}{r} \left[ 
\begin{picture}(7,2)
\put(0,-1){\includegraphics[width=7em]{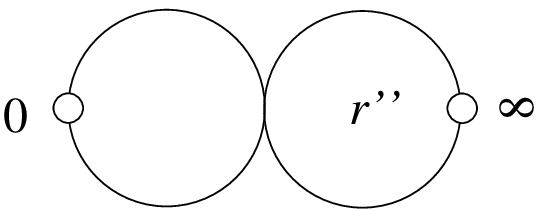}}
\end{picture}
\right], 
\end{equation}
where the sum is over all boundary divisors and $r''$ is the number of branch points on the component of~$\infty$.

\subsubsection{Computing $p^* \psi_s$}

Now we prove Theorems~\ref{Thm:geomlemma} and~\ref{Thm:geomlemmabis} using the preceding lemmas to express $p^*\psi_s$ in terms of boundary divisors. Assume for definiteness that $a_s > 0$. Then we have
$$
a_s p^* \psi_s = a_s \Psi_s - a_s(\Psi_s - p^* \psi_s) = 
q^* \psi_0 - a_s(\Psi_s - p^* \psi_s).
$$

Equations~(\ref{Eq:psibranch}) and~(\ref{Eq:psimarked}) give two alternative expressions for $q^* \psi_0$ while Lemma~\ref{Lem:Psiminuspsi} gives an expression for $a_s(\Psi_s - p^* \psi_s)$. All three expressions involve very similar summations over the same set of divisors, but with different coefficients.

\paragraph{Proof of Theorem~\ref{Thm:geomlemma}.}
We use Equation~(\ref{Eq:psibranch}) for $q^* \psi_0$. The coefficient of
$$
\frac{\prod_{i=1}^p k_i}{p!}
\left[\oM_{g_1;a_I,-k_1, \dots, -k_p}\right]^\virt 
\boxtimes 
\left[ \oM_{g_2; a_J, k_1, \dots, k_p} \right]^\virt
$$
in Eq.~(\ref{Eq:psibranch}) equals $r''/r$. Its coefficient in Lemma~\ref{Lem:Psiminuspsi} multiplied by $a_s$ equals~$1$ if $s \in J$ or $0$ if $s \in I$. Subtracting the second coefficient from the first one and using $r'+r''=r$ we get
$$
\frac{r''}{r} \quad \mbox{if} \quad s \in I,
$$
$$
-\frac{r'}{r} \quad \mbox{if} \quad s \in J.
$$
These are exactly the coefficients from Theorem~\ref{Thm:geomlemma}. \qed

\paragraph{Proof of Theorem~\ref{Thm:geomlemmabis}.}
We use Equation~(\ref{Eq:psimarked}) for $q^* \psi_0$. Denote by $l$ the index of the marked point with $a_l=0$ that appears in this equation. The coefficient of
$$
\frac{\prod_{i=1}^p k_i}{p!}
\left[\oM_{g_1;a_I,-k_1, \dots, -k_p}\right]^\virt 
\boxtimes 
\left[ \oM_{g_2; a_J, k_1, \dots, k_p} \right]^\virt
$$
in Eq.~(\ref{Eq:psimarked}) equals $1$ if $l \in J$ and $0$ otherwise. Its coefficient in Lemma~\ref{Lem:Psiminuspsi} multiplied by $a_s$ equals~$1$ if $s \in J$ and otherwise. Subtracting the second coefficient from the first we get
$$
1 \quad \mbox{if} \quad s \in I, l \in J,
$$
$$
-1\quad \mbox{if} \quad s \in J, l \in I,
$$
and $0$ otherwise. These are exactly the coefficients from Theorem~\ref{Thm:geomlemmabis}. \qed

Both theorems are proved.

\subsection{A digression on admissible coverings} 
Double ramification cycles have an alternative definition, using admissible coverings rather than relative stable maps (see, for instance,~\cite{Ion02}).
To distinguish the two notions, just for the length of this section, we will write $\DRa$ and $\DRs$. The goal of this section is to explain what would change in our results if we replaced $\DRs$ by $\DRa$. This section is not self-contained, since we don't introduce the admissible coverings here; it can be skipped in first reading.

\begin{example} \label{Ex:stabvsadm}
We have
$$
\DRa_1(a,\widetilde{-a}) = a^2 - 1 \in H^0(\oM_{1,1}),
$$
$$
\DRs_1(a,\widetilde{-a}) = a^2 \in H^0(\oM_{1,1}),
$$
where the tilde means that the the corresponding marked point is forgotten. 

Indeed, given an elliptic curve $(C,x)$ with one marked point, there exists $a^2$ points $y$ such that $x-y$ is an $a$-torsion point in the Jacobian of~$C$. The space of admissible coverings contains one point per $y \not= x$, that is, $a^2-1$ points. The space of rubber maps contains one additional point corresponding to $y=x$: it represents the map with a contracted elliptic component.
\end{example}

\begin{theorem} \label{Thm:stabvsadm}
The intersection numbers of a monomial $\psi_1^{d_1} \cdots \psi_n^{d_n}$ with $\DRa_g(a_1, \dots, a_n)$ and with $\DRs_g(a_1, \dots, a_n)$ coincide if none of the $a_i$'s vanishes. These intersection numbers may differ in presence of an $a_i = 0$.
\end{theorem}

\begin{corollary} \label{Cor:notpolynomial}
At least for some $g$ and~$n$ the class $\DRa_g(a_1, \dots, a_n)$ does not have a polynomial dependence on $a_1, \dots, a_n$.
\end{corollary}

\paragraph{Proof.} Our formulas show that the intersection number of a given monomial in $\psi$-classes with $\DRs_g(a_1, \dots, a_n)$ depends polynomially on $a_1, \dots, a_n$. The intersection number of the same monomial with $\DRa_g(a_1, \dots, a_n)$ has the same values for non-zero $a_i$'s, but different values if some of the $a_i$'s vanish. Therefore this intersection number cannot depend polynomially on $a_1, \dots, a_n$. Hence the class itself cannot depend polynomially on $a_1, \dots, a_n$. \qed

\begin{remark}
Ultimately it's Corollary~\ref{Cor:notpolynomial} that convinced us that $\DRs$-cycles must be preferred to $\DRa$-cycles.
\end{remark}

\paragraph{Proof of Theorem~\ref{Thm:stabvsadm}.}

\paragraph{A.} The first claim of the theorem is proved by checking that all the steps of our computation of $\DR_g(a_1, \dots, a_n) \psi_1^{d_1} \cdots \psi_n^{d_n}$ go through in the same way for $\DRa$ and $\DRs$ as long as $a_1, \dots, a_n$ do not vanish. 

Theorem~\ref{Thm:geomlemma} is the base of our computations. Its analogue for $\DRa$-cycles is well-known (see~\cite{Ion02}, Lemma~2.4 for a proof modulo some omitted terms; see~\cite{ShaZvo} Lemmas~3.2, 3.3, 3.6, 3.7 for a detailed proof is genus~1 that easily generalizes to higher genus). The proof is actually even simpler for $\DRa$-cycles because the space of admissible coverings has the expected dimension, so there is no virtual fundamental class involved. The combinatorial part of the computation only uses Theorem~\ref{Thm:geomlemma} as long as there are no vanishing $a_i$'s, therefore it works in the same way for both $\DRa$ and $\DRs$. 

The only part of the computation that does not generalize is the use of Okounkov and Pandharipande's computation in the case where there are $\psi$-classes only at the marked point with vanishing $a_i$'s. If there are no vanishing $a_i$'s this part is not needed.

\paragraph{B.} To prove the second claim of the theorem we will use the following example. Let
$$
\beta = \left[\includegraphics[scale=0.35,viewport=0 30 200 100]{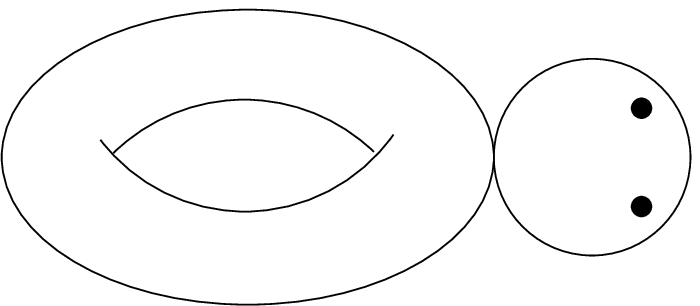}\right] ,
$$
$$
\gamma = \left[\includegraphics[scale=0.35,viewport=0 30 115 100]{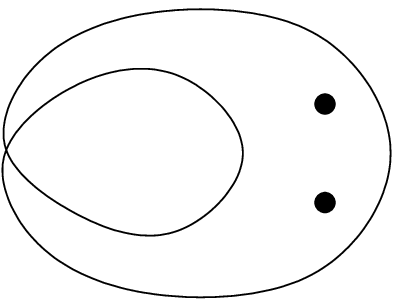}\right]
$$
be two cohomology classes in $H^2(\oM_{1,2})$.

\begin{proposition} We have
$$
\DRa_1(a,-a) = (a^2-1) \beta + \frac{a^2-1}{12} \gamma,
$$
$$
\DRs_1(a,-a) = a^2 \beta + \frac{a^2-1}{12} \gamma.
$$
\end{proposition}

\paragraph{Proof.} In this case the space of rubber maps has two irreducible components of the same dimension equal to the expected dimension. Its virtual fundamental class is the sum of the fundamental classes of the two components. The first component coincides with the space of admissible coverings. The second component is composed of maps with a contracted torus in the source curve. Thus the difference between $\DRa_1(a,-a)$ and $\DRs_1(a,-a)$ comes from the component with contracted tori, whose fundamental class projects to $\beta$. In other words, $\DRs(a,-a) - \DRa(a,-a) = \beta$, which is actually the only thing that is needed for the proof of the theorem.

The expression for $\DRs_1(a_1, \dots, a_n)$ is given by Hain~\cite{Hain11} in full generality, so our expression can be found as a particular case. Both formulas can also be proved using a lifting of the WDVV relation in the Losev-Manin space. \qed

\begin{corollary}
We have
$$
\DRa_1(a,-a,0) \psi_3^2 = (a^2-1)/12,
$$
$$
\DRs_1(a,-a,0) \psi_3^2 = (2a^2-1)/24.
$$
\end{corollary}

\paragraph{Proof.} The classes $\DR_1(a,-a,0)$ are obtained from $\DR_1(a,-a)$ by pull-backs under the forgetful map $\pi^* \colon \oM_{1,3} \to \oM_{1,2}$. It is straightforward to compute the intersection of the classes thus obtained with $\psi_3^2$. Note that the second equality is a particular case of Example~\ref{Ex:nequals3}. What matters for us is that
$$
\psi_3^2 (\DRs(a,-a,0) - \DRa(a,-a,0) = \psi_3^2 \pi^*\beta = \frac1{24} \ne 0.
$$
\qed

\bigskip

Thus we have found an example where the intersection numbers of the same monomial in $\psi$-classes with a $\DRa$-cycles and with a $\DRs$-cycle differ. This proves the second claim of the theorem. \qed

\section{Generating functions for one $\psi$-class}

In this section we prove Theorem~\ref{Thm:1} that evaluates the power of a $\psi$-class on a DR-cycle. We also generalize it to DR-cycles with forgotten points. From this we deduce Theorem~\ref{Thm:evidence}.

\subsection{Proof of Theorem~\ref{Thm:1}}

The proof of Theorem~\ref{Thm:1} splits into two very different cases: $a_s =0$ and $a_s \not= 0$. 

\subsubsection{The case $a_s =0$.} 

\begin{lemma} \label{Lem:psiequalspsi}
Let $p \colon \oM_{g; a_1, \dots, a_n} \to \oM_{g,n}$ be the forgetful map from the rubber space to the moduli space. Assume that $a_s = 0$. Then we have $p^* \psi_s = \psi_s$.
\end{lemma}

\paragraph{Proof.} Let $f \colon C \to \CP^1$ be a point of $\oM_{g;a_1, \dots, a_n}$. If the $s$th marked point lies on a component of the source curve contracted by~$f$ then this component is stable, because $f$ is stable. If it lies on a component that is not contracted by~$f$ then this component contains at least two more marked points: a pre-image of~$0$ and a pre-image of~$\infty$. Therefore it is also stable. Thus the $s$th marked point never lies on a component of the source curve contracted by the forgetful map. This allows us to identify the cotangent lines to the curve at the $s$th marked point before and after the forgetful map.\qed

Note that the statement of the lemma is completely wrong if $a_s \ne 0$.

\begin{lemma} \label{Lem:rubberequalsfixedimage}
Let $\mu, \nu$ be two partitions of the same integer~$d$. Consider the rubber space
$$
\oM^{\sim}_{g,p,\mu,\nu}(\CP^1,0,\infty)
$$
of relative stable maps to $\CP^1$ with $p$ marked points $x_1, \dots, x_p$.
Also consider the moduli space
$$
\oM^1_{g,p,\mu,\nu}(\CP^1,0,\infty)
$$
of relative stable maps to $\CP^1$ with $p$ marked points $x_1, \dots, x_p$ such that the image of $x_1$ is fixed to be $1 \in \CP^1$.
These two spaces are isomorphic to each other; their perfect obstruction theories and virtual fundamental classes coincide.
\end{lemma}

This is a well-known fact and the proof is a simple check.

The consequence of these two lemmas is that the intersection number
$$
\psi_s^{2g-3+n} \DR_g(a_1, \dots, a_n)
$$
is actually a Gromov-Witten invariant of $\CP^1$ relative to two points. Indeed, by Lemma~\ref{Lem:psiequalspsi}, instead of evaluating $\psi_s^{2g-3+n}$ we can evaluate it directly on the rubber space $\oM_{g; a_1, \dots, a_n}$. And, according to Lemma~\ref{Lem:rubberequalsfixedimage}, this is equivalent to finding a Gromov-Witten invariant of $\CP^1$ relative to two points. 

The Gromov-Witten invariants that we need were computed in~\cite{OkoPan06}, Eq.~(3.11) and \cite{Rossi}, Eq.~(10) if $a_s$ is the only vanishing marking; while the general case is covered by Proposition~2.5 of~\cite{OkoPan06b}. We do not have any new contribution to this computation. 

This proves the theorem for $a_s = 0$.

\subsubsection{Proof for $a_s \ne 0$.}

We proceed by induction on the number of branch points $r=2g-2+n_++n_-$. The base case is $r=1$, that is, $g=0$, $n_++n_-=3$. (Genus~1 is impossible, because $n_+ + n_- \geq 2$.) We have
$$
\psi_s^{n_0} \DR_0(a_1,a_2,a_3, \underbrace{0, \dots, 0}_{n_0}) = \int_{\oM_{0,n_0+3}} \psi_s^{n_0} = 1,
$$
which coincides with the constant term of the generating function in the theorem.

Recall the recursion from Corollary~\ref{Cor:onepsi}:
\begin{align}
&r a_s \; \psi_s^{r-1} \DR_g(a_1, \dots, a_n)
=  \label{formula4}\\
&- \frac12 \sum_{i,j \not= s}
(a_i+a_j) \; \psi_s^{r-2} \DR_g(a_1, \dots, \widehat{a_i}, \dots, \widehat{a_j}, \dots, a_n, a_i+a_j)  \label{formula5}\\
&-\frac12 \sum_{i \not= s} 
{\rm sign}(a_i) \sum_{\substack{b+c=a_i\\ b \cdot c>0}}
bc \; \psi_s^{r-2} \DR_{g-1}(a_1, \dots, \widehat{a_i}, \dots, a_n, b, c). \label{formula6}
\end{align}
By the induction assumption, the sum \eqref{formula5} is equal to
\begin{align*}
&-[z^{2g}]\frac{1}{ra_s}\sum_{i,j \not= s} \frac{\sinh\left(\frac{a_iz}{2}\right)\cosh\left(\frac{a_jz}{2}\right)}{\frac{z}{2}}\frac{\prod_{l\ne i,j,s}S(a_lz)}{S(z)}=\\
&=[z^{2g}]\frac{1}{ra_s}\sum_{i \not= s} (a_i+a_s)\cosh\left(\frac{a_iz}{2}\right)\frac{\prod_{j\ne i,s}S(a_jz)}{S(z)}.
\end{align*}
By the induction assumption, the sum \eqref{formula6} is equal to
\begin{align*}
&-[z^{2g-2}]\frac{1}{ra_s}\sum_{i \ne s} 
{\rm sign}(a_i) \sum_{\substack{b+c=a_i\\bc>0}}\frac{bc}{2}S(bz)S(cz)\frac{\prod_{l\ne i,s} S(a_lz)}{S(z)}=\\
&=-[z^{2g}]\frac{1}{ra_s} \sum_{i \ne s}
\left(a_i\cosh\left(\frac{a_iz}{2}\right)-\frac{\sinh\left(\frac{a_iz}{2}\right)
\cosh\left(\frac{z}{2}\right)}{\sinh\left(\frac{z}{2}\right)}\right)\frac{\prod_{j\ne i,s}S(a_jz)}{S(z)}.
\end{align*}
Thus, \eqref{formula4} is equal to
\begin{align*}
&[z^{2g}]\frac{1}{ra_s}
\sum_{i \ne s}
\left(a_s\cosh\left(\frac{a_iz}{2}\right)+
\frac{\sinh\left(\frac{a_iz}{2}\right)\cosh\left(\frac{z}{2}\right)}{\sinh\left(\frac{z}{2}\right)}\right)\frac{\prod_{j\ne i,s}S(a_jz)}{S(z)}=\\
&=[z^{2g}]\frac{1}{r}\left(z\frac{d}{dz}+n-2\right)\frac{\prod_{i \ne s} S(a_iz)}{S(z)}=\\
&=[z^{2g}]\frac{\prod_{i \ne s}S(a_iz)}{S(z)}.
\end{align*}
The theorem is proved. \qed

\begin{remark} \label{Rem:cutandjoin}
Our formula for $\psi_s^{r-1} \DR(a_1, \dots, a_n)$ coincides, up to a simple factor, with the formula for one-part double Hurwitz numbers found in~\cite{GoJaVa} (first equality of Theorem~3.1). This is due to the fact that the recursion relation of Corollary~\ref{Cor:onepsi} coincides with the cut-and-join equation for Hurwitz numbers. Note, however, that in the Hurwitz numbers theory the formula only holds for {\em one-part} numbers; in other words, all the numbers $a_1, \dots, a_n$ must be of the same sign except for $a_s$ that is of the opposite sign. If this condition is not satisfied then the cut-and-join equation fails. In our situation, however, the signs of the numbers $a_i$ do not matter. Thus we get the same generating function and the same cut-and-join equation, but their interpretations and their ranges of applicability are different.
\end{remark}

\subsection{Proof of Theorem~\ref{Thm:evidence}}

The statement of Theorem~\ref{Thm:evidence} involves a combination of intersection numbers of the form
$$
\psi_x^{m-K+p} \; 
DR_g\left((\widetilde{-1})^p, \framebox{p-K}, k_1, \dots, k_n\right),
$$
where the tilde means that the corresponding marked points are forgotten and the boxed marked point carries the $\psi$-class. So far we have obtained a formula evaluating a power of a $\psi$-class over a DR-cycle without forgotten points. So our first task is to generalize it to a DR-cycle with forgotten points.

\begin{proposition} \label{Prop:forgottenpoints}
Let $a_1, \dots, a_n$ and $b_1, \dots, b_p$ be two lists of integers such that $\sum a_i + \sum b_j = 0$. Let $s \in \{ 1, \dots, n \}$.
We have
$$
\psi_s^{2g-3+n+p} \DR_g(a_1, \dots, a_n; \; \tb_1, \dots, \tb_p) 
= [z^{2g}] \frac{\prod\limits_{i \ne s} S(a_i z) \prod\limits_{j=1}^p (S(b_j z) - 1)}
{S(z)}.
$$
\end{proposition}

\begin{example}
If $p>g$ the virtual dimension of the DR-space exceeds the dimension of $\oM_{g,n}$, so the image of the virtual fundamental class under the projection to $\oM_{g,n}$ vanishes. Thus the corresponding DR-cycle is equal to zero. The generating function has $p$ factors of the form $S(b_jz)-1$ that are divisible by $z^2$, therefore its coefficient of $z^{2g}$ also vanishes.
\end{example}

\begin{example}
Suppose that $b_p=0$. If $\pi \colon \oM_{g,n+1} \to \oM_{g,n}$ is the forgetful map, we have
$$
\DR(a_1, \dots, a_n; \; \tb_1, \dots, \tb_{p-1}, {\tilde 0}) =
\pi_* \pi^* \DR(a_1, \dots, a_n; \; \tb_1, \dots, \tb_{p-1}) = 0.
$$
As for the generating function, it contains a factor $S(0)-1 = 0$, so it vanishes identically.
\end{example}

\begin{example}
Let $p=g$ and $b_1 = \dots = b_g = 1$. Then we have
$$
\DR_g(a_1, \dots, a_n; {\tilde 1}, \dots, {\tilde 1}) = g! \in H^0(\oM_{g,n}).
$$
Indeed, it is known that the map from ${\rm Sym}^g C \to {\rm Pic}^g(C)$ has degree~1 for any smooth curve~$C$. Thus we have
$$
\psi_s^{3g-3+n} \DR_g(a_1, \dots, a_n; {\tilde 1}, \dots, {\tilde 1})
= g! \psi_s^{3g-3+n} \oM_{g,n} = \frac{g!}{g! \; 24^g} = \frac1{24^g}.
$$
The generating function contains $g$ factors 
$$
S(0)-1 = \frac{z^2}{24} + \cdots,
$$ 
thus 
$z^{2g}/24^g$ is its first non-zero term.

More generally, it is easy to show that 
$$
\DR_g(a_1, \dots, a_n; \tb_1, \dots, \tb_g) = g! \; b_1^2 \cdots b_g^2 \in H^0(\oM_{g,n})
$$
and therefore 
$$
\psi_s^{3g-3+n} \DR_g(a_1, \dots, a_n; \tb_1, \dots, \tb_g)
= \frac{b_1^2 \cdots b_g^2}{24^g}.
$$
\end{example}

\begin{lemma}
For any $j \in \{1, \dots, p \}$ we have 
$$
\psi_s^{2g-3+n+p} \; \DR_g(a_1, \dots, a_n; \; \tb_1, \dots, \tb_p)
$$
\begin{align*}
&= \psi_s^{2g-3+n+p} \; \DR_g(a_1, \dots, a_n, b_j; \;  \tb_1, \dots, \widehat{\tb_j}, \dots, \tb_p) 
\\
& - \psi_s^{2g-4+n+p} \; \DR_g(a_1, \dots, a_s+b_j, \dots, a_n; \;  \tb_1, \dots, \widehat{\tb_j}, \dots, \tb_p),
\end{align*}
where a hat means a skipped element.
\end{lemma}

\paragraph{Proof.} Let $\pi \colon \oM_{g,n+1} \to \oM_{g,n}$ be the forgetful map, where the $(n+1)$st marked point corresponds to the marking $b_j$. We have
$$
\pi^*\psi_s = \psi_s - D_{s,n+1},
$$
where the divisor $D_{s,n+1} \subset \oM_{g,n+1}$ encodes curves with a rational component containing exactly one node and two marked points with markings $s$ and~$n+1$. From this we deduce
$$
\psi_s^{2g-3+n+p} \; \DR_g(a_1, \dots, a_n; \; \tb_1, \dots, \tb_p)
$$
\begin{align*}
& = \pi^*\psi_s^{2g-3+n+p} \; \DR_g(a_1, \dots, a_n, b_j; \; \tb_1, \dots, \widehat{\tb_j}, \dots, \tb_p)
\\
& = (\psi_s - D_{s,n+1}) \pi^*\psi_s^{2g-4+n+p} \; \DR_g(a_1, \dots, a_n, b_j; \; \tb_1, \dots, \widehat{\tb_j}, \dots, \tb_p)
\\
& = \left( \psi_s^{2g-3+n+p} - D_{s,n+1} \pi^* \psi_s^{2g-4+n+p} \right)
\; \DR_g(a_1, \dots, a_n, b_j; \; \tb_1, \dots, \widehat{\tb_j}, \dots, \tb_p)
\\
&= \psi_s^{2g-3+n+p} \; \DR_g(a_1, \dots, a_n, b_j; \;  \tb_1, \dots, \widehat{\tb_j}, \dots, \tb_p) 
\\
& - \psi_s^{2g-4+n+p} \; \DR_g(a_1, \dots, a_s+b_j, \dots, a_n; \;  \tb_1, \dots, \widehat{\tb_j}, \dots, \tb_p).
\end{align*}
\qed

\paragraph{Proof of Proposition~\ref{Prop:forgottenpoints}.}
Applying the lemma successively for each $j \in \{ 1, \dots, p \}$ we get the following expression:
$$
\psi_s^{2g-3+n+p} \; \DR_g(a_1, \dots, a_n; \; \tb_1, \dots, \tb_p) 
$$
$$
= 
[z^{2g}] \sum_{J \subset \{1, \dots, p \}} (-1)^{p-|J|}
\frac{\prod_{i \ne s} S(a_i z) \prod_{j \in J} S(b_j z)}{S(z)}.
$$
This obviously factorizes into the expression given in the proposition. \qed

\paragraph{Proof of Theorem~\ref{Thm:evidence}.}
The statement of the theorem immediately follows from Proposition~\ref{Prop:forgottenpoints}. We have

\begin{multline*}
\sum_{p=0}^g
\frac{m!}{p! \, (K-p)!} \, \psi_x^{m-K+p} \;
\prod_{i=1}^n k_i\; 
\DR_g\left(\framebox{K-p}, -k_1, \dots, -k_n; \; {\tilde 1}^p\right)
\hspace{-20em}\\
= m! \frac{\prod_{i=1}^n k_i}{K!} 
\sum_{p=0}^g
\binom{K}{p} \psi_x^{m-K+p} \;
\DR_g\left(\framebox{K-p}, -k_1, \dots, -k_n; \; {\tilde 1}^p\right)
\\
=m! \frac{\prod_{i=1}^n k_i}{K!} 
[z^{2g}] \frac{\prod_{i=1}^n S(k_i z)}{S(z)} \sum _{p=0}^g
\binom{K}{p} \left( S(z) - 1 \right)^p
\\
=m! \frac{\prod_{i=1}^n k_i}{K!} [z^{2g}] S(z)^{K-1} \prod_{i=1}^n S(k_i z).
\end{multline*}
This completes the proof of the theorem.
\qed

\section{Semi-infinite wedge formalism}\label{sectionInfiniteWedge}

In this section we sketch the theory of the semi-infinite wedge space  following \cite{OkoPan06},~\cite{Joh10} and~\cite{ShaSpiZvo}.

\subsection{The infinite wedge space}

Let $V$ be an infinite dimensional vector space with basis labelled by the half integers. Denote by $\underline{i}$ the basis vector labelled by $i$, so $V = \bigoplus_{i \in \Z + \frac{1}{2}} \underline{i}$.

\begin{definition}\label{Def:WedgeProduct} Let $c$ be an integer.
An {\em infinite wedge product of charge $c$} is a formal expression 
\begin{equation}\label{eq:wedgeProduct}
\underline{i_1} \wedge \underline{i_2} \wedge \cdots
\end{equation}
such that the sequence of half-integers $i_1, i_2, i_3, \dots$ differs from the sequence $c-1/2, c-3/2, c-5/2, \dots$ in only a finite number of places.

The infinite wedge product is anti-commutative under finite permutations of factors.
\end{definition}

\begin{definition}
The {\em charged infinite wedge space} is the span of all infinite wedge products. The {\em infinite wedge space} is its zero charge subspace, that is, the span of all zero charge infinite wedge products.
\end{definition}

The infinite wedge space is spanned by the vectors
$$
v_\lambda = \underline{\lambda_1 - 1/2} \wedge \underline{\lambda_2 - 3/2} \wedge \underline{\lambda_3 - 5/2} \wedge \dots,
$$
where $(\lambda_1 \geq \lambda_2 \geq \dots \geq 0 \geq 0 \geq \dots)$ is a partition of any non-negative integer.

On both spaces we introduce the inner product $(,)$ for which the vectors~(\ref{eq:wedgeProduct}) are orthonormal.

\begin{definition}
The zero charge vector $v_\emptyset = \underline{-\frac{1}{2}} \wedge \underline{-\frac{3}{2}} \wedge \cdots$ is denoted by $|0 \rangle$ and is called the {\em vacuum vector}. Its dual $\langle 0 |$ with respect to the inner product is called the {\em covacuum} vector.
\end{definition}

\begin{definition}
If $\mathcal{P}$ is an operator on the infinite wedge space, we define its  {\em vacuum expectation value} as
$\langle \mathcal{P}\rangle^\circ = \langle 0 |\mathcal{P}|0\rangle$.
\end{definition}

\subsection{The operators}

\begin{definition} Let $k$ be any half integer. Then the operator $\Psi_k$ is defined by
$\Psi_k \colon (\underline{i_1} \wedge \underline{i_2} \wedge \cdots) \ \mapsto \ (\underline{k} \wedge \underline{i_1} \wedge \underline{i_2} \wedge \cdots)$. This operator acts on the charged infinite wedge space and increases the charge by $1$.

The operator $\Psi_k^*$ is defined to be the adjoint of the operator $\Psi_k$ with respect to the inner product.
\end{definition}

The action of $\Psi_k^*$ on basis vectors can be described as follows: it looks for a factor $\underline{k}$ in the wedge product, anti-commutes it to the beginning of the product, and erases it. If the factor $\underline{k}$ is not found the operator returns~$0$. Thus it decreases the charge by~1.

\begin{definition}
The normally ordered products of $\Psi$-operators are defined in the following way
\begin{equation}
{:}\Psi_i \Psi_j^*{:} \ = \begin{cases}\Psi_i \Psi_j^*, & \text{ if } j > 0 \\
-\Psi_j^* \Psi_i & \text{ if } j < 0\ .\end{cases} 
\end{equation}
\end{definition}

Note that the two expressions are equal unless $i=j$. Also note that the operator~${:}\Psi_i \Psi_j^*{:}$ does not change the charge of an infinite wedge product, and can thus be viewed as an operator on the infinite wedge space.

\begin{definition}
Let $n \not=0$ be a non-zero integer. We define the operator $\alpha_n$ on the infinite wedge space by
$$
\alpha_n = \sum_{k \in \Z + \frac{1}{2}} {:} \Psi_{k-n} \Psi^*_k {:} \; .
$$
\end{definition}

The operator $\alpha_n$ attempts to increase every factor of an infinite wedge product by~$n$ and returns the sum of successful attempts.

\begin{remark}
Let $E_{ij}$ for $i, j \in \Z + \frac{1}{2}$ denote the standard basis of matrix units of $\mathfrak{gl}(\infty) = \mathfrak{gl}(V)$. Then the assignment $E_{ij} \mapsto\ {:}\Psi_i \Psi_j^*{:}$ defines a projective representation of the Lie algebra $\mathfrak{gl}(V)$ on the infinite wedge space.
\end{remark}

%
%

These operators relate the infinite wedge space to the representation theory of the symmetric group and to double Hurwitz numbers via the following properties.

\begin{proposition} \label{Prop:character}
Let $\nu = (\nu_1, \dots, \nu_k)$ be a partition of~$N$. Then we have
$$
\prod_{i=1}^k \alpha_{-\nu_i} | 0 \rangle = \sum_{\lambda}
\chi_\lambda(\nu) v_\lambda,
\qquad
\langle 0 | \prod_{i=1}^k \alpha_{\nu_i} | v_\lambda \rangle = \chi_\lambda(\nu),
$$
where $\chi_\lambda(\nu)$ is the character of the conjugacy class $\nu$ in the representation~$\lambda$ and is set to be $0$ if $\nu$ and $\lambda$ are partitions of two different integers.
\end{proposition}

This proposition is a reformulation of the Murnaghan-Nakayama rule~\cite{James}.

\begin{proposition} \label{Prop:id}
We also have
$$
{\rm id} = \sum_{k \geq 1} \frac1{k!} \sum_{\nu_1, \dots, \nu_k} 
\frac1{\prod \nu_i} |\prod \alpha_{\nu_i} | 0 \rangle \langle 0 | \prod \alpha_{-\nu_i} |,
$$
where id the identity operator in the infinite wedge space.
\end{proposition}

This proposition follows from the orthogonality of characters.

The operator $\alpha_n$ will be used to ``create'' a branch point of order~$n$ over~0 (if $n>0$) or a branch point of order $-n$ over $\infty$ (if $n<0$). Insertions of the identity operator in the form of Proposition~\ref{Prop:id} can be used to describe a pinching of the target sphere separating $0$ from $\infty$.


\begin{notation}\label{notationZeta}
We denote by $\zeta(z)$ the function $e^{z/2} - e^{-z/2}$.
\end{notation}

\begin{definition}
Let $n \in \Z$ be any integer. We define the operator $\E_n(z)$ depending on a formal variable~$z$ by
$$
\E_n(z) = \sum_{k \in \Z + \frac{1}{2}} e^{z(k - \frac{n}{2})} \, {:} \Psi_{k-n} \Psi^*_k {:} \, + \frac{\delta_{n,0}}{\zeta(z)}.
$$
\end{definition}

Note that $\E_n(0) = \alpha_n$ for $n \not=0$.

\begin{proposition}\label{Prop:Ecommut}
We have
$$
[\E_k(w), \E_l(z)] = \zeta(kz-lw) \E_{k+l}(z+w);
$$
in particular,
$$
[\E_k(0) , \E_l(z)] = \zeta(kz) \E_{k+l}(z)
$$
and, taking a limit as $z \to 0$,
$$
[\E_k(0) , \E_l(0)] = k  \delta_{k+l,0}.
$$
\end{proposition}
Most of the time these commutation relations will be the only properties of the operators $\E$ that we are going to use.

\begin{remark} \label{Rem:Ereal}
Note that the definition of $\E_n(z)$ makes perfect sense for any {\em real} number~$n$ and the commutation relations above remain valid. Most of the time we will only use integer values of~$n$, but occasionally we will have to consider a limit as $n \to 0$ of an expression containing~$\E_n$. 
\end{remark}

\subsection{Commutation patterns}

Now we briefly describe an algorithm to compute a vacuum expectation value 
$$
\langle \E_{a_1}(z_1) \cdots \E_{a_n}(z_n) \rangle^\circ
$$
of a product of $\E$ operators, where $\sum a_i = 0$. The algorithm is described in detail in~\cite{Joh10} and~\cite{ShaSpiZvo}.

\begin{definition}\label{Def:Energy}
We say that the operator $\mathcal{E}_n(z)$ has positive (resp., negative, zero) energy if the integer~$n$ is positive (resp., negative, zero). 
\end{definition}

\begin{proposition}
We have $\E_n(z)|0\rangle = 0 $ for $n>0$, and $\langle 0|\E_n(z) = 0$ for $n<0$.
\end{proposition}


\begin{example}
Here is how we use the above proposition to compute vacuum expectation values:
$$
\begin{array}{ccccc}
\langle \E_3(z) \E_{-3}(w) \rangle^\circ & = & 
\langle \E_{-3}(z) \E_3(w)  \rangle^\circ & + &
\zeta(3z+3w) \langle \E_0(z+w)\rangle^\circ\\
\\
& = & 0 & + & \displaystyle \frac{\zeta(3(z+w)}{\zeta(z+w)}.
\end{array}
$$
Here the first equality is the commutation of $\E_3(z)$ and $\E_{-3}(w)$, while the second equality follows from the proposition and the definition of $\E_0$.
\end{example}

In general, the algorithm consists in ``moving'' operators of negative energy to the left. At each step we replace two successive operators by the sum of their commutator (the {\em commutator term}) and their product in the inverse order (the {\em passing term}). Both terms of this sum are still vacuum expectation values of products of $\E$-operators, multiplied by some product of $\zeta$-functions. As soon as an operator of negative energy reaches the leftmost position in the product we can cancel the corresponding term in the sum.

It is easy to see that this algorithm always terminates in a finite number of steps leaving us with a sum of terms of the form
$$
\zeta(\cdot) \cdots \zeta(\cdot) \langle \E_0(\cdot) \cdots \E_0(\cdot) \rangle^\circ.
$$
Since $\E_0(z) v_\emptyset = 1/\zeta(z) v_\emptyset$, the last expression is immediately computable.

\begin{definition}
Every term in the final sum together with the sequence of steps of the algorithm leading to it is called a {\em commutation pattern}.
\end{definition}

\begin{definition}
A commutation pattern is called {\em connected} if there is only one $\E_0$ operator in the brackets of the corresponding final term.
\end{definition}

\begin{definition}
A {\em connected vacuum expectation value} is the sum of outcomes of all connected commutation patterns.
\end{definition}

\begin{notation}
A connected expectation value will be denoted by~$\langle \cdot \rangle$, while a possibly disconnected expectation value retains the notation~$\langle \cdot \rangle^\circ$. 
\end{notation}

\begin{proposition}
The connected expectation value is well-defined, that is, it does not depend on the order in which the commutations of the algorithm are performed.
\end{proposition}

\paragraph{Sketch of a proof.} Every $\E_0$ factor in the final bracket of a commutation pattern arises from a combination of several factors $\E_{a_i} (x_i)$ for $i \in I$ such that $\sum_{i \in I} a_i = 0$. In particular, if there are no non-trivial subsets $I \subset \{ 1, \dots, n \}$ such that $\sum_{i \in I} a_i=0$ the connected expectation value coincides with the possibly disconnected expectation value. In general, the connected expectation value can be obtained from the possibly disconnected expectation values using the exclusion-inclusion formula. Thus it can be defined unambiguously without invoking our algorithm. \qed

\section{Proof of Theorem~\ref{Thm:2}}

In this section we use the infinite wedge formalism to prove Theorem~\ref{Thm:2}. The proof is by induction, and the base case is Theorem~\ref{Thm:1}. For that we first need to restate both theorems in terms of the infinite wedge formalism. 

In the rest of this section $g$ will always be used to denote the genus of some $\DR$-cycle which is intersected with some $\psi$-classes. It is determined by the dimension constraint that this intersection should be a number. 

\begin{proposition}\label{prop:1psiclass} Let $a_1, \dots, a_n$ be a list of real numbers such that $\sum a_i=0$. Denote $J = \{1, \ldots n\}\backslash \{s\}$ and define 
$$
J_+ = \{\ j \in J \ : \ a_j \geq 0\ \}
$$ 
and $J_- = J \backslash J_+$. Then
\begin{equation}
\label{eqn:1psiclass}
[z^{2g}] \frac{\prod_{j \in J} S(a_j z)}{S(z)}
= [x^{2g-3+n}]
\left\langle \prod_{j \in J_+} \frac{\E_{a_j}(0)}{a_j} \E_{a_s}(x) \prod_{j \in J_-}\frac{\E_{a_j}(0)}{-a_j} \right\rangle.
\end{equation}
\end{proposition}

\begin{remark}
When~$a_i = 0$ for some~$i$, we interpret the right-hand side of the formula as the limit for $a_i \to 0$ (cf Remark~\ref{Rem:Ereal}). Note that the apparent singularity coming from the factor~$\frac{1}{a_i}$ will be cancelled with the zero coming from~$\zeta(a_i x)$ when we compute the vacuum expectation value.
\end{remark}

\begin{lemma} Let $a_1, \dots, a_n$ be a list of non-zero real numbers with vanishing sum. Assume that it is split into a disjoint union of three sets
$$
\{1 , \dots, n \} = \{ s \} \sqcup J_+ \sqcup J_-.
$$
Then the vacuum expectation value 
$$
\left\langle \prod_{j \in J_+} \E_{a_j}(0) \E_{a_s}(x) \prod_{j \in J_-} \E_{a_j}(0) \right\rangle
$$
vanishes unless all the elements of $J_+$ are positive and all the elements of $J_-$ are negative.
\end{lemma}

\paragraph{Proof.} Assume, for instance, that an element $a_j = a$ of $J_-$ is positive. We will apply our algorithm by moving this operator to the right. Since $[\E_a(0), \E_b(0)] = a \delta_{0,a+b}$, we see that the commutator term either vanishes or contains an $\E_0$ factor that is prohibited in a connected expectation value. Thus to compute the connected expectation value we only need to take the passing term. In other words we can just move the operator $\E_a(0)$ to the right-most position and we see that the expectation value vanishes, because $\E_a(0) v_\emptyset = 0$.
\qed

\paragraph{Proof of Proposition~\ref{prop:1psiclass}.} We can assume that the $a_i$'s do not vanish since the case of vanishing $a_i$'s is obtained by passing to the limit.

To compute the expectation value we apply our algorithm by commuting $\E_{a_s}(x)$ with its right and left neighbours in an arbitrary order. At every step the contribution of the passing term vanishes by the lemma. Using the commutation formulas we get 
\begin{multline*}
[x^{2g - 3 + n}] \left\langle  \E_0(x) \right\rangle \prod_{j \in J_+} \frac{\zeta(a_j x)}{a_j} \prod_{j \in J_-}\frac{\zeta(-a_j x)}{-a_j} \\
=  [x^{2g+n-2}] \frac{x^{n-1}\prod_{j \in J} S(a_j x)}{\zeta(x)} = [x^{2g}]\frac{\prod_{i \neq s}S(a_i x)}{S(x)} ,
\end{multline*}
where we used $S(0) = 1$ and $S(-a_jx) = S(a_j x)$. \qed

\bigskip

Proposition~\ref{prop:1psiclass} allows us to restate Theorem~\ref{Thm:1} in terms of the infinite wedge formalism:

\begin{corollary} \label{Cor:thm1bis}
We have
$$
\psi_s^{2g-3+n} \DR_g(a_1, \dots, a_n) = [x^{2g-3+n}]
\left\langle \prod_{j \in J_+} \frac{\E_{a_j}(0)}{a_j} \E_{a_s}(x) \prod_{j \in J_-}\frac{\E_{a_j}(0)}{-a_j} \right\rangle.
$$
\end{corollary}

On the other hand, Theorem~\ref{Thm:2} is equivalent to the following:

\begin{theorem}\label{Thm:general}
Let $n$ be a positive integer, and let $a_1, \ldots, a_n$ be integers such that $\sum a_i = 0$ and $d_1, \dots, d_n$ non-negative integers such that $\sum d_i = 2g-3+n$. Then we have
\begin{multline}\label{eqn:dr-iw}
\psi_1^{d_1} \cdots \psi_n^{d_n} \DR_g(a_1, \ldots, a_n)  \\
= [x_1^{d_1} \cdots x_n^{d_n}] \sum_{\sigma \in S'_n} \frac{\left\langle \left. \left[\cdots\left[\E_{a'_1}(x'_1), \E_{a'_2}(x'_2)\right], \ldots \right], \E_{a'_n}(x'_n)\right] \right\rangle}{x'_1(a'_2 - \frac{a'_1x'_2}{x'_1})\cdots(a'_n - \frac{a'_{n-1}x'_{n}}{x'_{n-1}})} ,
\end{multline}
where $S'_n$ is the group of permutations~$\sigma$ of the set $\{1, \ldots, n\}$ with $\sigma(1) = 1$. As before, we define $x'_i := x_{\sigma(i)}$ and $a'_i := a_{\sigma(i)}$. 
\end{theorem}
The fact that Theorem~\ref{Thm:general} is equivalent to Theorem~\ref{Thm:2} follows immediately by repeated application of the commutation relation of Proposition~\ref{Prop:Ecommut}.

\begin{notation}
In the following, if a sequence of integers~$a_1, \ldots, a_n$ and a corresponding sequence of formal variables~$x_1, \ldots, x_n$ have been defined, we will often use the following abbreviation for the sake of clarity: 
$$
\E_j := \E_{a_j}(x_j) .
$$
\end{notation}

\begin{definition}
Let $I \subset \{1, \dots, n \}$.
Given a list of integers $a_1, \dots, a_n$ and a list of variables $x_1, \dots, x_n$, let $P$ be a polynomial in operators $\E_i$, $i \in I$, whose coefficients are rational functions in $x_i$ and $a_i$. Let $t \not\in I$.

For any $i \in I$ we define $\mathcal{O}^t_i P$ to be the result of the substitution
$$
\E_i \mapsto \frac{1}{\frac{a_i x_t}{x_i} - a_t} [\E_i, \E_t]. 
$$
Furthermore, we define $\mathcal{O}^t P = \sum_{i \in I} \mathcal{O}^t_i P$. 
\end{definition}

\begin{definition}\label{defG}
Let $n$ be some positive integer, let $a_1, \ldots, a_n$ be some sequence of integers, and let $x_1, \ldots, x_n$ be a sequence of formal variables. Then we define
$$
G^t(a_1, \ldots, a_n; x_1, \ldots, x_n) = \sum_{\sigma \in S'_t} \frac{\left. \left[\cdots\left[\E_{\sigma_1}, \E_{\sigma_2}\right], \ldots \right], \E_{\sigma_t}\right]}{x_{\sigma_1} (\frac{a_{\sigma_1} x_{\sigma_2}}{x_{\sigma_1}} - a_{\sigma_2}) \cdots (\frac{a_{\sigma_{t-1}} x_{\sigma_t}}{x_{\sigma_{t-1}}} - a_{\sigma_t})} .
$$
\end{definition}

To prove Theorem~\ref{Thm:general}, we will need the following lemmas. 

\begin{lemma}\label{operatorLemma}
For any positive integers~$t \leq n$, and for all $a_1, \ldots, a_n \in \Z$, we have
\begin{equation}\label{operEq}
G^t(a_1, \ldots , a_n; x_1, \ldots, x_n) = \mcO^t \cdots \mcO^2 \frac{1}{x_1}\E_1 .
\end{equation}
Note that the empty product of operators that appears on the right-hand side in the case $t=1$ should be interpreted as the identity operator, as usual. 
\end{lemma}

\paragraph{Proof.} We prove that the coefficients of $x_1^{d_1} \cdots x_n^{d_n}$ are equal on both sides of the equation for any non-negative integers~$d_1, \ldots, d_n$. The lemma is clearly true when~$t = 1$. We proceed by induction on~$t$. 

Denote by $F^t$ the right-hand side of equation~\eqref{operEq}:
$$
F^t := \mcO^t \cdots \mcO^2 \frac{1}{x_1}\E_1 \\
$$
Now assume that $F^t$ and~$G^t$ are equal for some~$t$, and are related by just a series of applications of the Jacobi identity. Defining
$$
\tilde{G}^{t+1}_i = \mcO^{t+1}_i G^t
$$
it follows that $\mcO^{t+1}_i F^t = \tilde{G}^{t+1}_i$, again related by a series of applications of the Jacobi identity. We complete the proof by showing that $G^{t+1}$ is equal to $\tilde{G}^{t+1} := \sum_{i=0}^{t} \tilde{G}^{t+1}_i$, and this equality can be given just by application of the Jacobi identity. 

The terms of $G^{t+1}$ are of the form 
\begin{equation}\label{termOfG}
\frac{[\cdots [\E_{\sigma_1}, \E_{\sigma_2}], \ldots ],\E_{\sigma_i}], \E_{t+1}], \ldots],\E_{\sigma_t}]}{x_{\sigma_1} (\frac{a_{\sigma_i}x_{t+1}}{x_{\sigma_i}} - a_{t+1})(\frac{a_{t+1}x_{\sigma_{i+1}}}{x_{t+1}} - a_{\sigma_{i+1}})\prod_j (\frac{a_{\sigma_j} x_{\sigma_{j+1}}}{x_{\sigma_j}} - a_{\sigma_{j+1}})} ,
\end{equation}
where~$\sigma$ is some permutation appearing in the sum in~$G^t$, and where~$j$ runs from $1$ to~$t-1$, skipping~$i$. 

First we look at the case where $0<i<t$. A term with the iterated commutator appearing in~\eqref{termOfG} arises in $\tilde{G}^{t+1}$ in precisely two ways:
\begin{enumerate}
\item In $\tilde{G}^{t+1}_i$, as the first term of the Jacobi identity applied to
$$
[\cdots [\E_{\sigma_1}, \E_{\sigma_2}], \ldots ],\E_{\sigma_{i-1}}], [\E_{\sigma_i}, \E_{\sigma_{t+1}}],\E_{\sigma_{i+1}}], \ldots],\E_{\sigma_t}]
$$
\item In $\tilde{G}^{t+1}_{i+1}$, as the second term of the Jacobi identity applied to
$$
[\cdots [\E_{\sigma_1}, \E_{\sigma_2}], \ldots ],\E_{\sigma_{i}}], [\E_{\sigma_i+1}, \E_{\sigma_{t+1}}],\E_{\sigma_{i+2}}], \ldots],\E_{\sigma_t}].
$$
\end{enumerate}
Taking into account the coefficients of these two contributions, we get
\begin{multline*}
\frac{1}{\prod_j (\frac{a_{\sigma_j} x_{\sigma_{j+1}}}{x_{\sigma_j}} - a_{\sigma_{j+1}})} \left(\frac{1}{(\frac{a_{\sigma_i} x_{\sigma_{i+1}}}{x_{\sigma_i}} - a_{\sigma_{i+1}}) (\frac{a_{\sigma_i} x_{t+1}}{x_{\sigma_i}} - a_{t+1})} \right.\\
\left. - \frac{1}{(\frac{a_{\sigma_i} x_{\sigma_{i+1}}}{x_{\sigma_i}} - a_{\sigma_{i+1}}) (\frac{a_{\sigma_{i+1}} x_{t+1}}{x_{\sigma_{i+1}}} - a_{t+1})}\right) \\
= \frac{1}{\prod_j (\frac{a_{\sigma_j} x_{\sigma_{j+1}}}{x_{\sigma_j}} - a_{\sigma_{j+1}})} \frac{1}{\frac{a_{\sigma_i}x_{t+1}}{x_{\sigma_i}} - a_{t+1}}\frac{1}{\frac{a_{t+1}x_{\sigma_{i+1}}}{x_{t+1}} - a_{\sigma_{i+1}}} ,
\end{multline*}
(where~$j$ runs over the same set as above) which is precisely the coefficient of the iterated commutator appearing in~\eqref{termOfG}. 

In this way, we use all the terms of~$\tilde{G}^{t+1}$ except three to get the terms of the form~\eqref{termOfG} in~$G^{t+1}$ with~$0<i<t$. The three terms we did not yet use are~$\tilde{G}^{t+1}_0$, the second term of the Jacobi identity applied to $\tilde{G}^{t+1}_1$ and the first term of the Jacobi identity applied to~$\tilde{G}^{t+1}_t$. The only remaining terms in~$G^{t+1}$ are those of the form~\eqref{termOfG} with~$i=0$ and with~$i=t$. By a similar argument as the one above, the first is easily seen to be equal to the sum of $\tilde{G}^{t+1}_0$ and the second term of the Jacobi identity applied to $\tilde{G}^{t+1}_1$, whereas it is 
immediate that the second is equal to the first term of the Jacobi identity applied to~$\tilde{G}^{t+1}_t$. This completes the induction. \qed

\begin{corollary}\label{cor:order}
Let $\bar{a} = (a_1, \ldots, a_n)$ be an ordered set of positive integers. Denote $\bar{x} = (x_1, \ldots, x_n)$. Then the expression 
$$
F(\bar{x};\bar{a}) := \mcO^n \cdots \mcO^2 \frac{1}{x_1} \E_{a_1}(x_1)
$$
is symmetric with respect to the action of~$S_n$ on~$\{1, \ldots, n\}$. 
\end{corollary}

\paragraph{Proof.} The group $S_n$ is generated by $S'_n$ and the transposition~$(1,2)$. The fact that $F(\bar{x},\bar{a})$ is symmetric with respect to the action of~$S'_n$ follows immediately because the expression on the right-hand side of~\eqref{operEq} is symmetric with respect to this action. The invariance under the transposition~$(1,2)$ is shown as follows:
$$
\mcO^2 \frac{1}{x_1}\E_1 = \frac{[\E_1, \E_2]}{a_1 x_2 - a_2 x_1} = \frac{[\E_2,\E_1]}{x_2(\frac{a_2 x_1}{x_2} - a_1)} = \mcO^1 \frac{1}{x_2}\E_2. 
$$
\qed

\begin{lemma}\label{lem:nopsi}
Let~$n$ be any positive integer, and let~$a_1, \ldots, a_n$ be positive integers. For any subset~$I \subset \{2, \ldots, n\}$ we have
$$
[\prod_{i \in I^\mathrm{c}} x_i^0] \left\langle \mcO^n \cdots \mcO^2 \frac{1}{x_1} \E_1\right\rangle = \left\langle \prod_{i \in I^\mathrm{c}, a_i \geq 0} \frac{\E_i}{a_i} \left( \prod_{t \in I} \mcO^t \frac{1}{x_1} \E_1\right)\prod_{j \in I^\mathrm{c}, a_j <0} \frac{\E_j}{-a_j} \right\rangle ,
$$
where~$I^\mathrm{c}$ denotes the complement of~$I \subset \{2, \ldots, n\}$. 
\end{lemma}

\paragraph{Proof.}
Let $k = |I|$. By Corollary~\ref{cor:order}, we can assume that $I = \{1, \cdots, k\}$. 

First note that
$$
[\prod_{i \in I^\mathrm{c}} x_i^0] \mcO^n \cdots \mcO^2 \frac{1}{x_1} \E_1 = \frac{\tilde{\mcO}^n}{a_n}  \cdots \frac{\tilde{\mcO}^{k+1}}{a_{k+1}} \mcO^k \cdots \mcO^2 \frac{1}{x_1} \E_1
$$
where $\tilde{\mcO}^t$ is defined by $\tilde{\mcO}^t := \mathrm{ad}_{\E_{a_t}(0)}$. This follows immediately from the fact that we take the coefficient of $x_{k+1}^0 \cdots x_n^0$ and the expansion of the factor~$\frac{1}{\frac{a_i x_t}{x_i} - a_t}$ appearing in the operator~$\mcO^t_j$.

Whenever $a_t = 0$ for some $t > k$, the operator $\frac{\tilde{\mcO}_t}{a_t}$ just acts as multiplication by~$\frac{\E_t}{a_t}$ from the left; since we take the connected vacuum expectation value, it will be forced to involved in a commutator with a negative energy operator from the left at some point. It is easy to see that the effect of this commutator will be the same as the effect from the one coming from the definition of~$\tilde{\mcO}_t$. 

For $t > k$ with $a_t \neq 0$, we use the standard Lie-theory fact that for any $t > k$
\begin{multline*}
\tilde{\mcO}^t \cdots \tilde{\mcO}^{k+1} \mcO^k \cdots \mcO^2 \frac{1}{x_1} \E_1 = \E_{a_t}(0) \tilde{\mcO}^{t-1} \cdots \tilde{\mcO}^{k+1} \mcO^k \cdots \mcO^2 \frac{1}{x_1} \E_1 \\
- \tilde{\mcO}^{t-1} \cdots \tilde{\mcO}^{k+1} \mcO^k \cdots \mcO^2 (\frac{1}{x_1} \E_1) \E_{a_t}(0) . 
\end{multline*}
Depending on the sign of~$a_t$, only one of these two terms will contribute when we take the vacuum expectation value. Iterating this procedure from $t = k + 1$ up to $t = n$ completes the proof of the lemma. \qed

\begin{remark}
Using Lemmas~\ref{operatorLemma} and~\ref{lem:nopsi}, it is easy to see that Corollary~\ref{Cor:thm1bis} is a special case of Theorem~\ref{Thm:general}, i.e. that Theorem~\ref{Thm:1} is a special case of Theorem~\ref{Thm:2}.
\end{remark}

\begin{lemma}\label{Cor:operatorMove}
Let $n$ be a positive integer, let $a_1, \ldots, a_n$ be a sequence of integers and let $x_1, \ldots, x_n$ be a sequence of formal variables. Denote $\bar{a} := (a_1, \ldots, a_n)$ and $\bar{x} = (x_1, \ldots, x_n)$. Then we have, for any $p, q$ with $1 \leq p < q \leq n$:
$$
G^n(\bar{a};\bar{x}) = \frac{x_p x_q}{a_p x_q - a_q x_p} \sum_{I,J} [G^{|I|}(a_I; x_I), G^{|J|}(a_J; x_J)]  
$$
where the sum is over all disjoint sets $I$ and~$J$ such that $I \cup J = \{1, \ldots, n\}$ and $p \in I$, $q \in J$, and $G$ is as defined in definition~\ref{defG}. 
\end{lemma}
\paragraph{Proof.} 
Let us introduce the following notation. Suppose $h_1,h_2,\ldots,h_n$ are operators. Let
$$
Q_n(h_1,\ldots,h_n;x_1,\ldots,x_n)=\sum_{\sigma \in S'_n} \frac{x_{\sigma_2}x_{\sigma_3}\cdots x_{\sigma_{n-1}}\left[ \left[\cdots\left[h_{\sigma_1}, h_{\sigma_2}\right], \ldots \right], h_{\sigma_n}\right]}{(x_{\sigma_1}-x_{\sigma_2}) \cdots (x_{\sigma_{n-1}}-x_{\sigma_n})}. 
$$
For any $\sigma\in S_n$, we have the symmetry 
$$
Q_n(h_{\sigma_1},\ldots,h_{\sigma_n};x_{\sigma_1},\ldots,x_{\sigma_n})=Q_n(h_1,\ldots,h_n;x_1,\ldots,x_n) .
$$
It can be proved in the same way as the symmetry of $G^n(\bar a;\bar x)$. We have
$$
G^n(\bar{a};\bar{x})=\frac{(-1)^{n-1}}{a_1a_2\ldots a_n}Q_n(\E_1,\ldots,\E_n;\frac{x_1}{a_1},\ldots,\frac{x_n}{a_n}).
$$
The lemma obviously follows from the formula
\begin{gather}\label{formula: tmp}
Q_n(h;z)=\frac{x_px_q}{x_p-x_q}\sum_{\substack{I\coprod J=\{1,\ldots,n\}\\p\in I, \; q\in J}}[Q_{|I|}(h_I;x_I),Q_{|J|}(h_J;x_J)].
\end{gather}
We prove \eqref{formula: tmp} by induction on $n$. The case $n=2$ is obvious. Suppose $n\ge 3$. Let us denote the set $\{1,2,\ldots,n\}$ by $\overline n$. We have
\begin{align}
&\frac{x_px_q}{x_p-x_q}\sum_{\substack{I\coprod J=\overline n\\p\in I\\q\in J}}[Q_{|I|}(h_I;x_I),Q_{|J|}(h_J;x_J)]=\notag\\
&=\sum_{i\ne p,q}\frac{x_p(x_i-x_q)}{(x_p-x_q)(x_p-x_i)}\times\notag\\
&\times\frac{x_i x_q}{x_i-x_q}\sum_{\substack{I\coprod J=\overline n\\i,p\in I\\q\in J}}[Q_{|I|-1}([h_p,h_i],h_{I\backslash \{i,p\}};x_i,x_{I\backslash\{i,p\}}),Q_{|J|}(h_J;x_J)]+\label{term 1}\\
&+\frac{x_q}{x_p-x_q}[h_p,Q_{n-1}(h_{\overline n\backslash\{p\}};x_{\overline n\backslash\{p\}})].\label{term 2}
\end{align} 
By the induction assumption the sum~\eqref{term 1} is equal to
\begin{align}
&\sum_{i\ne p,q}\frac{x_p(x_i-x_q)}{(x_p-x_q)(x_p-x_i)}Q_{n-1}([h_p,h_i],h_{\overline n\backslash\{i,p\}};x_i,x_{\overline n\backslash\{i,p\}})=\notag\\
&=\sum_{i\ne p,q}\frac{x_i}{x_p-x_i}Q_{n-1}([h_p,h_i],h_{\overline n\backslash\{i,p\}};x_i,x_{\overline n\backslash\{i,p\}})+\label{term 3}\\
&+\sum_{i\ne p,q}\frac{x_q}{x_q-x_p}Q_{n-1}([h_p,h_i],h_{\overline n\backslash\{i,p\}};x_i,x_{\overline n\backslash\{i,p\}})\label{term 4}
\end{align}
It is easy to see that \eqref{term 4} is equal to
\begin{align}
&\frac{x_q}{x_q-x_p}[h_p,Q_{n-1}(h_{\overline n\backslash\{p\}};x_{\overline n\backslash\{p\}})]+\label{term 5}\\
&+\frac{x_q}{x_p-x_q}Q_{n-1}([h_p,h_q],h_{\overline n\backslash\{p,q\}};x_q,x_{\overline n\backslash\{p,q\}}).\label{term 6}
\end{align}
Clearly, the sum of \eqref{term 3} and \eqref{term 6} is equal to $Q_n(h;x)$ and the sum of~\eqref{term 2} and~\eqref{term 5} is zero. The formula~\eqref{formula: tmp} is proved.   \qed

The proof of Theorem~\ref{Thm:2} is by induction, starting from two base cases. In the first case all $\psi$-classes are at points on the DR-cycle which are not mapped to zero or infinity. In the second case there is one point in the inverse image of zero or infinity were some non-zero power of a $\psi$-class appears, and all other $\psi$-classes again are at points which are not mapped to zero or infinity. The induction will then be completed using Corollary~\ref{Cor:movepsi}, which allows us to move $\psi$-classes between points on the boundary. We now prove the two base cases.

\begin{proposition}\label{prop:psiMiddle}
Let $n$ be some positive integers, and let $a_1, \ldots, a_n$ be integers. Let $d_1, \ldots, d_n$ be a set of non-negative integers such that $d_i$ is zero whenever $a_i \neq 0$. Then Theorem~\ref{Thm:general} holds. That is, if we denote $I = \{1 \leq i \leq n : d_i > 0 \}$, then
\begin{equation}\label{eq:middlePsi}
DR_g(a_1, \ldots, a_n) \prod_{i = 1}^n \psi_i^{d_i} = [\prod_{i=1}^n x_i^{d_i}] \la \mcO^n \cdots \mcO^2 \frac{1}{x_1} \E_1 \ra .
\end{equation}
\end{proposition}

\paragraph{Proof.} First note that by Lemma~\ref{operatorLemma}, equation~\eqref{eq:middlePsi} is indeed equivalent to the described special case of Theorem~\ref{Thm:general}. 


By Corollary~\ref{cor:order} we can assume that $1 \in I$ and that $i < j$ whenever $a_i = 0$ and $a_j \neq 0$. Let $t$ be the number of~$i$ for which $a_i = 0$. By Lemma~\ref{lem:nopsi}, the right-hand side of equation~\eqref{eq:middlePsi} is equal to
$$
[\prod_{i=1}^t x_i^{d_i}] \la \prod_{i, a_i \geq 0} \frac{\E_i}{a_i} \left( \mcO^t \cdots \mcO^2 \frac{1}{x_1}\E_1\right) \prod_{j, a_j <0} \frac{\E_j}{-a_j} \ra , 
$$
which shows that when~$I$ contains only one element, the statement of the proposition is a direct consequence of Corollary~\ref{Cor:thm1bis}. 

Furthermore, by~\cite{OkoPan06b}, Proposition~2.5;
$$
\DR_g(a_1, \ldots, a_n) \prod_{i \in I} \psi_i^{d_i} = \binom{d_1 + \cdots + d_n}{d_1, \ldots, d_n} \DR_g(a_1, \ldots, a_n) \psi_1^{d_1 + \cdots + d_n -n +1} .
$$
Using this and the expression above for the right-hand side of~\eqref{eq:middlePsi}, it only remains to show that 
$$
[x_1^{d_1} \cdots x_t^{d_t}] \mcO^t \cdots \mcO^2 \frac{1}{x_1}\E_1 = \binom{d_1 +\cdots + d_t}{d_1, \ldots, d_t} [x^{d_1 + \cdots d_t - t +1}] \frac{1}{x} \E_0(x) .
$$
It is a direct computation that this equation is equivalent to
$$
\mcO^t \cdots \mcO^2 \frac{1}{x_1} \E_1 = (x_1 + \cdots + x_t)^{t-2} \E_0(x_1 + \cdots + x_t) ,
$$
which is clearly true when $t = 1$, so we proceed by induction. Suppose that the equation above is true for all $t \leq l$ for some~$l \geq 1$.

Note that the action of~$\mcO^l_j$ on $\prod_i \mcO^i \E_1$ is given by the following actions.
\begin{itemize}
\item replace $x_j$ by $x_j + x_l$
\item replace $a_j$ by $a_j + a_l$
\item multiply the result by 
$\frac{x_j \zeta(a_jx_l - a_lx_j)}{a_jx_l - a_lx_j} = x_j S(a_j x_l - a_l x_j)$, which is equal to~$x_j$ when $a_j$ and~$a_l$ tend to~$0$ (which is the case for us since $i \in I$ implies $a_i = 0$).
\end{itemize}
Thus, we have, by the induction hypothesis
\begin{multline*}
\mcO^{l+1} \cdots \mcO^2 \E_{a_1}(x_1) = \mcO^{l+1} x_1 (x_1 + \cdots + x_l)^{l-2}\E_0(x_1 + \cdots + x_l) \\
= (x_1 + x_{l+1})(x_1 + \cdots + x_{l+1})^{l-2}\E_0(x_1 + \cdots + x_{l+1}) x_1 \\
+ \sum_{j=2}^l x_1(x_1 +\cdots + x_{l+1})^{l-2}\E_0(x_1 + \cdots x_{l+1}) x_j \\
= x_1(x_1 + \cdots + x_{l+1})^{l-1} \E_0(x_1 + \cdots + x_{l+1}) .
\end{multline*}
\qed

\begin{proposition}\label{prop:1psiMix}
Let $n$ be some positive integer, and let $a_1, \ldots, a_n$ be some sequence of integers with~$a_1 \neq 0$. Let $K$ be a subset of $\{2, \ldots, n\}$, and suppose that $a_i$ is zero for all~$i \in K$. Then for all sequences of integers~$d_1, \ldots, d_n$ with~$d_i \neq 0$ if and only if $i \in K \cup \{1\}$ we have
$$
\DR_g(a_1, \ldots, a_n) \prod_{i = 1}^n \psi_i^{d_i} 
= [\prod_{i = 1}^n x_i^{d_i}]\la \mcO^n \cdots \mcO^2 \frac{1}{x_1}\E_{a_1}(x_1) \ra. 
$$
where $g$ is determined by the usual formula.
\end{proposition}

\paragraph{Proof.} Let $k := |K|$. By Corollary~\ref{cor:order} we can reorder the set $\{1, \ldots, n\}$ in such a way that the elements of~$K$ correspond to~$2,\ldots, k$. By Lemma~\ref{lem:nopsi}, the statement of the proposition is equivalent to
\begin{multline}\label{prrt}
\sum_{d_1, \ldots, d_k} \DR_g(a_1, \ldots, a_n) \prod_{i=1}^k (\psi_ix_i)^{d_i}\\ 
=  \la \left(\prod_{i > k, a_i \geq 0}\frac{\E_i}{a_i}\right) \left( \mcO^k \cdots \mcO^2 \frac{1}{x_1}\E_1\right) \left(\prod_{j>k, a_j<0}\frac{\E_j}{-a_j}\right) \ra. 
\end{multline}
We prove this statement using induction on~$k$. If~$k = 0$, the statement is a direct consequence of Corollary~\ref{Cor:thm1bis}. On the other hand, the right-hand side of~\eqref{prrt} is equal to
\begin{multline*}
\frac{1}{x_1(\frac{a_1x_2}{x_1} - a_2)} \la \left( \prod_{i > k, a_i \geq 0}\frac{\E_i}{a_i}\right) \left( \mcO^k \cdots \mcO^3 [\E_1,\E_2]\right) \left(\prod_{j>k, a_j<0}\frac{\E_j}{-a_j}\right) \ra \\
= \frac{1}{a_1 x_2 - a_2x_1} \sum_{S\subset\{2, \ldots, k\}} \\
\la \left(\prod_{i > k, a_i \geq 0}\frac{\E_i}{a_i}\right) \left( [\prod_{s \in S} \mcO^s\E_1, \prod_{t \in S^{\mathrm{c}}} \mcO^t \E_2]\right) \left(\prod_{j>k, a_j<0}\frac{\E_j}{-a_j}\right) \ra\\
= \frac{1}{a_1 x_2 - a_2x_1} \sum \left\{\frac{k_1 \cdots k_p}{k!}\vphantom{\la \prod_{i \in I}\E_i \left( \prod_{s \in S} \mcO^s\E_1\right) \prod_{j\in J} \E_j \prod_{q>0} \E_{-c_q} \ra} \right.\\
\la \prod_{i \in I}\frac{\E_i}{a_i} \left( \prod_{s \in S} \mcO^s\E_1\right) \prod_{j\in J} \frac{\E_j}{-a_j} \prod_{q=1}^p \frac{\E_{-k_q}(0)}{k_q} \ra \cdot \phantom{abacadabra}\\
\phantom{abacadabra}\cdot \la \prod_{q=1}^p \frac{\E_{k_q}(0)}{k_q} \prod_{i \in I^{\mathrm{c}}}\frac{\E_i}{a_i} \left( \prod_{s \in S^{\mathrm{c}}} \mcO^s\E_2\right) \prod_{j\in J^{\mathrm{c}}} \frac{\E_j}{-a_j}  \ra \\
- \la \prod_{i \in I}\frac{\E_i}{a_i} \left( \prod_{s \in S} \mcO^s\E_2\right) \prod_{j\in J} \frac{\E_j}{-a_j} \prod_{q=1}^p \frac{\E_{-k_q}(0)}{k_q} \ra \cdot \phantom{abacadabra} \\ 
\phantom{abacadabra} \cdot \left. \la \prod_{q=1}^p \frac{\E_{k_q}(0)}{k_q} \prod_{i \in I^{\mathrm{c}}}\frac{\E_i}{a_i} \left( \prod_{s \in S^{\mathrm{c}}} \mcO^s\E_1\right) \prod_{j\in J^{\mathrm{c}}} \frac{\E_j}{-a_j}  \ra \right\}, 
\end{multline*}
where the last sum is over all subsets~$S$ as in the line above, all subsets $I \subset \{i ; i > k, a_i \geq 0\}$ and $J \subset \{j ; j > k, a_j < 0 \}$, all positive integers~$t$ and all sequences of integers $k_1, \ldots, k_p$. If we take the coefficient of $x_1^{d_1} \cdots x_n^{d_n}$ in this expression, using the induction hypothesis and Proposition~\ref{prop:psiMiddle}, we get
\begin{multline*}
\sum \frac{1}{a_1} \frac{k_1 \cdots k_p}{p!}
\big( \DR_{g_1}(a_1, \bar{a}_1, k_1, \ldots, k_p) \boxtimes \DR_{g_2}(a_2, \bar{a}_2, k_1, \ldots, k_p)  \\
-  \DR_{g_1}(a_2, \bar{a}_1, k_1, \ldots, k_p) \boxtimes \DR_{g_2}(a_1, \bar{a}_2, k_1, \ldots, k_p) \big) \psi_1^{d_1 - 1} \psi_2^{d_2} \cdots \psi_n^{d_n} 
\end{multline*}
where we use $\bar{a}_1$ as shorthand for the set of variables $\{a_i\}_{i \in I \cup J \cup S}$ and $\bar{a}_2$ as shorthand for $\{a_i\}_{i \in I^\mathrm{c} \cup J^\mathrm{c} \cup S^\mathrm{c}}$, and where the sum is over the same range as in the previous equation (note that the genera $g_1$ and~$g_2$ are determined by dimensional constraints). By Theorem~\ref{Thm:geomlemmabis}, this combination is precisely equal to the intersection of $\psi$-classes and a DR-cycle as in the proposition. \qed. 

\paragraph{Proof of Theorem~\ref{Thm:general}.} 
The left-hand side of equation~\eqref{eqn:dr-iw} is completely determined by the case where $a_i$ is zero whenever~$d_i$ is not zero, the case there is precisely one $i$ with $a_i \neq 0$ and $d_i \neq 0$, and Corollary~\ref{Cor:movepsi}. 

The right-hand side of the equation is completely determined by the same two cases, and Lemma~\ref{Cor:operatorMove}. By Proposition~\ref{prop:psiMiddle} the left- and right-hand side are equal in the first case, and by Proposition~\ref{prop:1psiMix} they are equal in the second case. 

Furthermore, intersecting the equation of Corollary~\ref{Cor:movepsi} with a monomial in $\psi$-classes of total degree~$2g - 4 + n$, we see that the application of Corollary~\ref{Cor:movepsi} and Corollary~\ref{Cor:operatorMove} lead to equivalent operations on the left- and right-hand sides of equation~\eqref{eqn:dr-iw}. \qed

\vspace{4 mm}
A.~B.: Korteweg-de~Vries Institute for Mathematics, University of Amsterdam, P.~O.~Box 94248, 1090 GE Amsterdam, The Netherlands 
\emph{and} 
Department of Mathematics, Moscow State University, Leninskie gory, 119992 GSP-2 Moscow, Russia

\emph{E-mail address:} a.y.buryak@uva.nl

\vspace{4 mm} 
S.~S. and L.~S.: Korteweg-de~Vries Institute for Mathematics, University of Amsterdam, P.~O.~Box 94248, 1090 GE Amsterdam, The Netherlands

\emph{E-mail address:} s.shadrin@uva.nl, l.spitz@uva.nl

\vspace{4 mm}

D.~Z: Institut Math\'{e}matique de Jussieu, CNRS and UPMC, 4 place Jussieu, 75013 Paris, France

\emph{E-mail address:} dimitri.zvonkine@gmail.com

\end{document}